\documentclass[11pt]{article}
\usepackage{amsmath}
\usepackage{color}
\usepackage[colorlinks=true]{hyperref}
\usepackage{url}
\def\uphi{\varphi}
\def\rg{{\hbox{\tiny\rm ,rg}}}
\def\lf{{\hbox{\tiny\rm ,lf}}}
\def\r{{\hbox{\tiny\rm ,r}}}
\def\l{{\hbox{\tiny\rm ,l}}}
\def\up{\hbox{\tt right}}
\def\dw{\hbox{\tt left}}

\def\I{{\cal I}}

\def\cX{{\cal X}}
\def\cN{{\cal N}}

\newcommand{\half}{ \mbox{\small$\frac{1}{2}$}}

\usepackage{amsfonts}
\usepackage{epsfig}

\def\Risk{{\hbox{\rm Risk}}}

\def\cO{{\cal O}}

\usepackage{amssymb}

\usepackage{graphicx}

\def\bR{{\mathbf{R}}}
\def\bZ{{\mathbf{Z}}}
\def\Opt{{\hbox{\rm Opt}}}

\def\cL{{\cal L}}

\def\cT{{\cal T}}
\def\cE{{\cal E}}
\def\cP{{\cal P}}

\def\cM{{\cal M}}
\def\cX{{\cal X}}

\def\cN{{\cal N}}

\def\cF{{\cal F}}
\definecolor{MyDarkBlue}{rgb}{0,0.08,0.45}
\definecolor{MyViolet}{rgb}{0.45,0.08,0.95}
\definecolor{MyBrown}{rgb}{0.45,0.08,0}

\def\bE{{\mathbf{E}}}

\def\Prob{\hbox{\rm Prob}}

\def\Argmin{\mathop{\hbox{\rm Argmin}}}
\newtheorem{lemma}{Lemma}[section]

\newtheorem{proposition}{Proposition}[section]

\newtheorem{example}{Example}[section]

\def\qed{$\Box$}

\def\e{{\rm e}}

\newcommand{\aic}[2]{{\color{MyDarkBlue}~#2}}

\newcommand{\be}{\begin{eqnarray}}
\newcommand{\ee}[1]{\label{#1}\end{eqnarray}}

\newcommand{\ese}{\end{eqnarray*}}
\newcommand{\bse}{\begin{eqnarray*}}
\newcommand{\rf}[1]{~(\ref{#1})}
\newcommand{\wh}[1]{{\widehat{#1}}}
\definecolor{MyDarkBlue}{rgb}{0,0.08,0.45}
\definecolor{MyViolet}{rgb}{0.45,0.08,0.95}
\definecolor{MyBrown}{rgb}{0.45,0.08,0}

\title{Near-Optimal Recovery of Linear and $N$-Convex Functions on Unions of Convex Sets}
\author{
Anatoli Juditsky
\thanks{LJK, Universit\'e Grenoble Alpes, 700 Avenue Centrale 38401 Domaine Universitaire
de Saint-Martin-d'H\`{e}res, France,
{\tt anatoli.juditsky@univ-grenoble-alpes.fr}}
\and Arkadi Nemirovski
\thanks{Georgia Institute
 of Technology, Atlanta, Georgia
30332, USA, {\tt nemirovs@isye.gatech.edu}\newline
The first author was supported by the LabEx PERSYVAL-Lab (ANR-11-LABX-0025) and the PGMO grant 2016-2032H. Research of
the second author was supported by NSF grant   CCF-1523768.}}
\date{}
\begin{document}
\maketitle
\begin{abstract}
In this paper we build provably near-optimal, in the minimax sense, estimates of linear forms and, more generally, ``$N$-convex functionals'' (an example being the maximum of several fractional-linear functions) of unknown ``signal'' from indirect noisy observations, the signal assumed to belong to the union of finitely many given convex compact sets. Our main assumption is that the observation scheme in question is {\sl good} in the sense of \cite{GJN2015}, the simplest example being the Gaussian scheme where the observation is the sum of linear image of the signal and the standard Gaussian noise. The proposed estimates, same as upper bounds on their worst-case risks, stem from solutions to explicit convex optimization problems, making the estimates ``computation-friendly.''
\end{abstract}

\section{Introduction}
The simplest version of the problem considered in this paper is as follows. Given access to $K$ independent observations
\begin{equation}\label{eq0}
\omega_t=Ax + \sigma\xi_t,\,1\leq t\leq K\quad [A\in\bR^{m\times n},\xi_t\sim\cN(0,I_m)]
\end{equation}
of ``signal'' $x$ known to belong to the union $X=\bigcup_{i=1}^IX_i$ of convex compact sets $X_i\subset\bR^n$, we want to recover $f(x)$, where $f$ is either linear, or, more generally, {\sl $N$-convex}. Here $N$-convexity means that $f:\cX\to\bR$ is a continuous function on a convex compact domain $\cX\supset X$ such that for every $a\in\bR$, each of the two level sets $\{x\in\cX: f(x)\geq a\}$ and $\{x\in\cX: f(x)\leq a\}$  can be represented as the union of at most $N$ convex compact sets\footnote{Immediate examples are affine-fractional functions $f(x)=(a^Tx+a)/(b^Tx+b)$ with denominators positive on $\cX$, in particular, affine functions ($N=1$), and piecewise linear functions like $
\max[a^Tx+a,\min[b^Tx+b,c^Tx+c]]$ ($N=3$). A less trivial example is conditional quantile of a discrete distribution $(N=2)$, see Section \ref{sect:setting}.}. Our principal contribution is an estimation routine which is provably near-optimal in the minimax sense. Our construction is not restricted to the Gaussian observation scheme (\ref{eq0}) and deals with {\sl good observation schemes}\footnote{Our main results can be easily extended to the more general case of {\em simple families} -- families of distributions specified in terms of upper bounds on their moment-generating functions, see \cite{PartI,PartII} for details. Restricting the framework to the case of good observation schemes is aimed at  streamlining the presentation.} (o.s.'s),
as defined in \cite{GJN2015}; aside of  the  Gaussian o.s., important examples are
\begin{itemize}
\item {\sl Poisson o.s.}, where $\omega_t$ are independent across $t$ identically distributed vectors with independent across $i\leq m$ entries $[\omega_t]_i\sim\hbox{\rm Poisson}(a_i^Tx)$, and
\item {\sl Discrete o.s.}, where $\omega_t$ are independent across $t$ realizations of discrete random variable taking values $1,...,m$ with probabilities affinely parameterized by $x$.
\end{itemize}
  The problem of (near-)optimal recovery of {\sl linear} function $f(x)$ on a convex compact set or a finite union of convex sets $X$  {has received much attention in the statistical literature (see, e.g., \cite{IKhas84,Donoho1987,donoho1991geom2,donoho1991geom3,donlow,Don95,cailow,cailow2004,cailow2005,JN2009}).} {In particular,} D. Donoho proved, see \cite{Don95}, that in the case of Gaussian observation scheme (\ref{eq0}) and convex and compact $X$, the worst-case, over $x\in X$, risk of the minimax optimal {\sl affine in observations estimate} is within factor 1.2 of the actual
   minimax risk.\footnote{Here risks are the mean square ones, see  \cite{Don95} for details.} Later,
in \cite{JN2009}, this near-optimality result was extended to other good observation schemes.  In \cite{cailow2004,cailow2005} the minimax affine estimator was used as  ``working horse'' to build the near-optimal estimator  of a linear functional over a finite union $X$ of convex compact sets in the  Gaussian observation scheme.
As compared to the existing results, our contribution here is twofold. First, we pass from Gaussian o.s. to essentially more general good o.s.'s, extending in this respect the results of  \cite{cailow2004,cailow2005}. Second, we  relax the requirement of affinity of the function to be recovered to $N$-convexity of the function.\par
It should be stressed that the actual ``common denominator'' of the cited contributions and of the present work is the ``operational nature'' of the results, as opposed to typical results of non-parametric statistics which can be considered as descriptive. The traditional results present near-optimal estimates and their risks  in a ``closed analytical form,'' the toll being
severe restrictions on the families $X$ of signals and  observation schemes. For instance, in the case of (\ref{eq0}) such ``conventional''  results would impose strong and restrictive assumptions on the interconnection between the geometries of $X$ and  $A$. In contrast, the approach we advocate here, same as that of, e.g., \cite{Don95,JN2009}, allows for quite general, modulo convexity, signal sets $X_i$, for arbitrary matrices $A$ in the case of (\ref{eq0}), etc., and the proposed estimators and their risks are yielded by efficient computation rather than being given in a closed analytical form. All we know in advance is that those risks are nearly as low as they can be under the circumstances.

  \par
  The main body of the paper is organized as follows.  Section \ref{sect:simpleschemes} contains preliminaries,  originating from \cite{JN2009,GJN2015}, on good o.s.'s. In Section \ref{sect:linear} we deal with recovery of linear functions on the unions of convex sets. Finally, recovery of $N$-convex functions is the subject of Section \ref{sect:recNconvex}.
  It is worth to mention that the construction of near-optimal estimator used in Section \ref{sect:recNconvex} is completely different from that employed in \cite{Don95,cailow,cailow2004,cailow2005,JN2009} and is closely related to the binary search estimator from \cite{Donoho1987,donoho1991geom2} dealing with what can be seen as continuous analogue of discrete o.s..
  \footnote{In the hindsight, it is interesting to note that the authors of \cite{donoho1991geom2} believed their ``... estimator not intended to be implemented on a computer...'' They considered their construction as purely theoretical and finally oriented their analysis in the ``traditional'' way, by imposing assumptions allowing to end up with explicit convergence rates in some specific situations.}
 \par
  Some technical proofs are relegated to Appendix.
\section{Preliminaries: good observation schemes}\label{sect:simpleschemes}
The estimates to be developed in this paper heavily exploit the notion of a good observation scheme introduced in \cite{GJN2015}. To make the presentation self-contained we start with explaining this notion here.
\subsection{Good observation schemes: definitions}
Formally, a {\sl good observation scheme} (o.s.) is a collection $\cO=\left((\Omega,P),\{p_\mu(\cdot):\mu\in\cM\},\cF\right)$, where
\begin{itemize}
\item $(\Omega,P)$ is an {\sl observation space}: $\Omega$ is a Polish (complete metric separable) space, and $P$ is a $\sigma$-finite $\sigma$-additive Borel reference measure on $\Omega$, such that $\Omega$ is the support of $P$;
\item $\{p_\mu(\cdot):\mu\in\cM\}$ is a parametric family of probability densities, specifically, $\cM$ is a convex relatively open set in some $\bR^M$, and for $\mu\in \cM$, $p_\mu(\cdot)$ is a probability density, taken w.r.t. $P$, on $\Omega$. We assume that the function $p_\nu(\omega)$ is positive and continuous in $(\mu,\omega)\in\cM\times\Omega$;
\item $\cF$ is a finite-dimensional linear subspace in the space of continuous functions on $\Omega$. We assume that $\cF$ contains constants and all functions of the form $\ln(p_\mu(\cdot)/p_\nu(\cdot))$, $\mu,\nu\in \cM$, and that the function
    \begin{equation}\label{PhicalO}
    \Phi_{\cO}(\phi;\mu)=\ln\left(\int_\Omega \e^{\phi(\omega)} p_\mu(\omega)P(d\omega)\right)
    \end{equation}
    is  real-valued on $\cF\times\cM$ and is {\sl concave} in $\mu\in \cM$; note that this function is automatically convex in $\phi\in\cF$. From real-valuedness, convexity-concavity and the fact that both $\cF$ and $\cM$ are convex and relatively open, it follows that $\Phi$ is continuous on $\cF\times\cM$.
\end{itemize}
\subsection{Examples of good observation schemes}
As shown in \cite{GJN2015} (and can be immediately verified), the following o.s.'s are good:
\begin{enumerate}
\item {\sl Gaussian o.s.}, where $P$ is the Lebesgue measure on $\Omega=\bR^d$, $\cM=\bR^d$, $p_\mu(\omega)$ is the density of the Gaussian distribution $\cN(\mu,I_d)$ (mean $\mu$, unit covariance), and $\cF$ is the family of affine functions on $\bR^d$. Gaussian o.s. with $\mu$ linearly parameterized by signal $x$ underlying observations, see (\ref{eq0}), is the standard observation model in signal processing;
\item {\sl Poisson o.s.}, where $P$ is the counting measure on the nonnegative integer $d$-dimensional lattice $\Omega=\bZ^d_+$,  $\cM=\bR^d_{++}=\{\mu=[\mu_1;...;\mu_d]>0\}$, $p_{\mu}$ is the density, taken w.r.t. $P$, of random $d$-dimensional vector with independent $\hbox{\rm Poisson}(\mu_i)$ entries, $i=1,..,d$, and $\cF$ is the family of all affine functions on $\Omega$. Poisson o.s. with $\mu$ affinely parameterized by signal $x$ underlying observation is the standard observation model in {\sl Poisson imaging}, including Positron Emission Tomography \cite{SheppVardi1985}, Large Binocular Telescope \cite{Telescope2000,TelescopeOSEM2000}, and Nanoscale Fluorescent Microscopy, a.k.a. Poisson Biophotonics
\cite{hell1994breaking,hell2003toward,betzig2006imaging,hell2008microscopy,hess2006ultra};
\item {\sl Discrete o.s.}, where $P$ is the counting measure on the finite set $\Omega=\{1,2,...,d\}$, $\cM$ is the set of positive $d$-dimensional probabilistic vectors $\mu=[\mu_1;...;\mu_d]$, $p_\mu(\omega)=\mu_\omega$, $\omega\in\Omega$, is the density, taken w.r.t. $P$, of a probability distribution $\mu$ on $\Omega$, and $\cF=\bR^d$ is the space of all real-valued functions on   $\Omega$;
\item {\sl Direct product of good o.s.'s}. Given $K$ good o.s.'s $\cO_t=((\Omega_t,P_t),\{p_{t,\mu}:\mu\in\cM_r\},\cF_t)$, $t=1,...,K$, we can build from them a new (direct product) o.s. $\cO_1\times....\times\cO_K$  with  observation space $\Omega_1\times...\times\Omega_K$, reference measure $P_1\times...\times P_K$, family of probability densities $\{p_\mu(\omega_1,...,\omega_K)=\prod_{t=1}^Kp_{t,\mu_t}(\omega_t):\mu=[\mu_1;...;\mu_K]\in\cM_1\times...\times\cM_K\}$, and $\cF=\{\phi(\omega_1,...,\omega_K)=\sum_{t=1}^K\phi_t(\omega_t):\phi_t\in\cF_t,t\leq K\}$. In other words, the direct product of o.s.'s $\cO_t$ is the observation scheme in which we observe collections $\omega^K=(\omega_1,...,\omega_K)$ with independent across $t$ components $\omega_t$ yielded by o.s.'s $\cO_t$.
    \par
    When all factors $\cO_t$, $t=1,...,K$, are identical to each other, we can reduce the direct product $\cO_1\times...\times\cO_K$ to its ``diagonal,''
    referred to as {\sl $K$-th power $\cO^K$}, or {\sl stationary $K$-repeated version}, of $\cO=\cO_1=...=\cO_K$. Just as in the direct product case, the observation space and reference measure in $\cO^K$ are $\Omega^K=\underbrace{\Omega\times...\times\Omega}_{K}$ and $P^K=\underbrace{P\times...\times P}_{K}$, the family of densities is $\{p^K_\mu(\omega^K)=\prod_{t=1}^Kp_\mu(\omega_t):\mu\in\cM\}$, and the family $\cF$ is $\{\phi^{(K)}(\omega_1,...,\omega_K)=\sum_{t=1}^K\phi(\omega_t):\phi\in\cF\}$. Informally, $\cO^K$ is the observation scheme we arrive at when passing from a single observation drawn from a distribution $p_\mu$, $\mu\in\cM$,  to $K$ independent observations drawn from the same distribution $p_\mu$.\par
    It is immediately seen that direct product of good o.s.'s, same as power of good o.s., are themselves good o.s.
\end{enumerate}
\section{Recovering linear forms on unions of convex sets}\label{sect:linear}
Our objective now is to extend the results of \cite{JN2009} to the situation where $X$ is finite union of convex sets. At the same time, the results of this section can be seen as an extension to more general observation schemes of the constructions of \cite{cailow2004,cailow2005}.
\subsection{The problem}\label{sect:linear:prob}
Let $\cO=\left((\Omega,P),\{p_\mu(\cdot):\mu\in\cM\},\cF\right)$ be a good o.s.. The problem we are interested in this section is as follows:
\begin{quote}
We are given a positive integer $K$ and $I$ nonempty convex compact sets $X_j\subset\bR^n$, along with affine mappings $A_j(\cdot):\bR^n\to \bR^M$ such that $A_j(x)\in\cM$ whenever $x\in X_j$, $1\leq j\leq I$. In addition, we are given a linear function $g^Tx$ on $\bR^n$.
\par
Given random observation
$$\omega^K=(\omega_1,...,\omega_K)
$$
with $\omega_k$ drawn, independently across $k$, from  $p_{A_j(x)}$ with $j\leq I$ and $x\in X_j$, we want to recover $g^Tx$. It should be stressed that both $j$ and $x$ underlying our observation are unknown to us.
\end{quote}
Given reliability tolerance $\epsilon\in(0,1)$, we quantify the performance of a candidate estimate -- a Borel function $\widehat{g}(\cdot):\Omega^K\to\bR$ -- by the worst case, over $j$ and $x$, width of $(1-\epsilon)$-confidence interval. Specifically, we say that $\widehat{g}(\cdot)$ is {\em $(\rho,\epsilon)$-reliable}, if
$$
\forall (j\leq I,x\in X_j): \Prob_{\omega^K\sim p^K_{A_j(x)}}\{|\widehat{g}(\omega^K)-g^Tx|>\rho\}\leq\epsilon.
$$
We define $\epsilon$-risk of the estimate as the smallest $\rho$ such that $\widehat{g}$ is $(\rho,\epsilon)$-reliable:
$$
\Risk_\epsilon[\widehat{g}]=\inf\left\{\rho: \widehat{g}\hbox{\ is $(\rho,\epsilon)$-reliable}\right\}.
$$

\subsection{The estimate}\label{sec:unionest}
Following \cite{JN2009}, we introduce  parameters $\alpha>0$ and $\phi\in \cF$, and associate with a pair $(i,j)$, $1\leq i,j\leq I$, the functions
$$
\begin{array}{rcl}
\Phi_{ij}(\alpha,\phi;x,y)&=&
\half K\alpha\left[\Phi_\cO(\phi/\alpha;A_i(x))+\Phi_\cO(-\phi/\alpha;A_j(y))\right]+\half g^T[y-x]+\alpha\ln(2I/\epsilon):\\
&&
\quad\quad\quad\quad\quad\quad\quad\quad\quad\quad\quad\quad\quad\quad\quad{\{\alpha>0,\phi\in\cF\}\times[X_i\times X_j]\to\bR,}\\
\Psi_{ij}(\alpha,\phi)&=&\max_{x\in X_i,y\in X_j}\Phi_{ij}(\alpha,\phi;x,y)=\half\left[\Psi_{i,+}(\alpha,\phi)+\Psi_{j,-}(\alpha,\phi)\right]:\{\alpha>0\}\times\cF\to\bR,\\
\end{array}
$$
where
\bse
\Psi_{\ell,+}(\beta,\psi)&=&\max_{x\in X_\ell}\left[K\beta\Phi_\cO(\psi/\beta;A_\ell(x))-g^Tx+\beta\ln(2I/\epsilon)\right]:\{\beta>0,\psi\in\cF\}\to\bR,\\
\Psi_{\ell,-}(\beta,\psi)&=&\max_{x\in X_\ell}\left[K\beta\Phi_\cO(-\psi/\beta;A_\ell(x))+g^Tx+\beta\ln(2I/\epsilon)\right]:\{\beta>0,\psi\in\cF\}\to\bR
\ese
and $\Phi_\cO$ is given by (\ref{PhicalO}).
Note that the function $\alpha\Phi_\cO(\phi/\alpha;A_i(x))$ is obtained from continuous convex-concave function $\Phi_\cO(\cdot,\cdot)$ by projective transformation in the convex argument, and affine substitution in the concave argument, so that the former function is convex-concave and continuous on the domain $\{\alpha>0,\phi\in\cF\}\times X_i$. By similar argument, the function $\alpha\Phi_{\cO}(-\phi/\alpha;A_j(y))$ is convex-concave and continuous on the domain $\{\alpha>0,\phi\in\cF\}\times X_j$. These observations combine with compactness of $X_i$ and $X_j$ to imply that $\Psi_{ij}(\alpha,\phi)$ is real-valued continuous convex function on the domain
$$\cF^+=\{\alpha>0\}\times\cF.$$
Observe that functions $\Psi_{ii}(\alpha,\phi)$ are 
positive on $\cF^+$. Indeed, for any $\bar{x}\in X_i$, when setting $\mu=A_i(\bar{x})$, we have
$$
\begin{array}{l}
\Psi_{ii}(\alpha,\phi)\geq \Phi_{ii}(\alpha,\phi;\bar{x},\bar{x}) = {\alpha\over 2}\left[K[\Phi_\cO(\phi/\alpha;\mu)+\Phi_\cO(-\phi/\alpha;\mu)]+2\ln(2I/\epsilon)\right]\\
={\alpha\over 2}\bigg[K\ln\bigg(\underbrace{\left[\int\exp\{\phi(\omega)/\alpha\}p_\mu(\omega)P(d\omega)\right]\,
\left[\int\exp\{-\phi(\omega)/\alpha\}p_\mu(\omega)P(d\omega)\right]}_{\hbox{\tiny (by the Cauchy inequality)\ }\geq [\int\exp\{{1\over 2}\phi(\omega)/\alpha\}\exp\{-{1\over 2}\phi(\omega)/\alpha\}p_\mu(\omega)P(d\omega)]^2=1}\bigg)+2\ln(2I/\epsilon)\bigg]\\
\geq \alpha\ln(2I/\epsilon)>0.
\end{array}
$$
\par
Functions $\Psi_{ij}$ give rise to convex and feasible optimization problems
\begin{equation}\label{Optij}
\Opt_{ij}=\Opt_{ij}(K)=\min_{\alpha,\phi}\left\{\Psi_{ij}(\alpha,\phi):(\alpha,\phi)\in\cF^+\right\}.
\end{equation}
By construction, $\Opt_{ij}$ is either a real, or $-\infty$; by the observation above, $\Opt_{ii}$ is nonnegative.
Our estimate is as follows.
\begin{enumerate}
\item For $1\leq i,j\leq I$, we select a feasible solutions $\alpha_{ij},\phi_{ij}$ to problems (\ref{Optij}) (the smaller the values of the corresponding objectives, the better) and set
    \begin{equation}\label{gij}
\begin{array}{rcl}
\rho_{ij}&=&\Psi_{ij}(\alpha_{ij},\phi_{ij})={1\over 2}\left[\Psi_{i,+}(\alpha_{ij},\phi_{ij})+\Psi_{j,-}(\alpha_{ij},\phi_{ij})\right]\\
\varkappa_{ij}&=&{1\over 2}\left[\Psi_{j,-}(\alpha_{ij},\phi_{ij})-\Psi_{i,+}(\alpha_{ij},\phi_{ij})\right]\\
g_{ij}(\omega^K)&=&\sum_{k=1}^K\phi_{ij}(\omega_k)+\varkappa_{ij}
\end{array}
\end{equation}
\item Given observation $\omega^K$, we specify the estimate $\widehat{g}(\omega^K)$ as follows:
\begin{equation}\label{rcij}
\begin{array}{c}
\widehat{g}(\omega^K)={1\over 2}\left[\min_{i\leq I}r_i+\max_{j\leq I} c_j\right]\hbox{\ with\ }
r_i=\max_{j\leq I}g_{ij}(\omega^K),\,
c_j=\min_{i\leq I}g_{ij}(\omega^K)\\
\end{array}
\end{equation}
\end{enumerate}
%
\begin{proposition}\label{proplin} For $i\in \{1,...,I\}$, let $\rho_{i}=\max_{1\leq j\leq I} \max[\rho_{ij},\rho_{ji}]$, and let
\[
\rho=\max_{i}\rho_i=\max_{1\leq i,j\leq I}\rho_{ij}.
\]
Assume that the density , taken w.r.t. $P^K$, of the distribution of the $K$-repeated observation $\omega^K$ is $p^K_{A_\ell(x)}$ for some $\ell\leq I$ and $x\in X_\ell$. Then
\[
\Prob_{\omega^K\sim p^K_{A_\ell(x)}}\{|\widehat{g}(\omega^K)-g^Tx|>\rho_\ell\}\leq\epsilon.
\]
As a result, the $\epsilon$-risk of the estimate we have built satisfies
\begin{equation}\label{riskupperbound}
\Risk_\epsilon\left[\widehat{g}(\cdot)\right]\leq \rho.
\end{equation}
\end{proposition}
See Section \ref{pr:00} for the proof.
%
\par
Observe that  properly selecting $\phi_{ij}$ and $\alpha_{ij}$ we can make the upper bound $\rho$ on the $\epsilon$-risk of the above estimate arbitrarily close to
$$
\Opt(K)=\max_{1\leq i,j\leq I}\Opt_{ij}(K).
$$
We are about to show that the quantity $\Opt(K)$ ``nearly lower-bounds'' the minimax optimal $\epsilon$-risk
$$
\Risk^*_\epsilon(K)=\inf_{\widehat{g}(\cdot)} \Risk_\epsilon[\widehat{g}],
$$
where the infimum is taken over all $K$-observation Borel estimates. The precise statement is as follows:
\begin{proposition}\label{propnearoptlin} In the situation of this section, let $\epsilon\in(0,1/2)$ and  $\bar{K}$ be a positive integer. Then
for every integer $K$ satisfying
$$
K>{2\ln(2I/\epsilon)\over\ln({[4\epsilon(1-\epsilon)]^{-1}})}\bar{K}
$$
one has
\begin{equation}\label{onehasone}
\Opt(K)\leq \Risk^*_\epsilon(\bar{K}).
\end{equation}
In addition, in the special case where for every $i,j$ there exists $\bar x_{ij}\in X_i\cap X_j$ such that $A_i(\bar x_{ij})=A_j(\bar x_{ij})$ one has
\begin{equation}\label{onehastwo}
K\geq\bar{K}\Rightarrow \Opt(K)\leq {2\ln(2I/\epsilon)\over\ln({[4\epsilon(1-\epsilon)]^{-1}})}\Risk^*_\epsilon(\bar{K}).
\end{equation}
\end{proposition}
See Section \ref{pr:01} for the proof.

\subsection{Illustration} We illustrate our construction by applying it to the simplest possible example in which the observation scheme is Gaussian and $X_i=\{x_i\}$ are singletons in $\bR^n$, $i=1,...,I$. Setting $y_i=A_i(x_i)\in \bR^m$, the observation components $\omega_k$, $1\leq k\leq K$, stemming from $(i,x_i)$, are drawn independently of each other from the normal distribution $\cN(y_i,I_m)$. Recall that in the Gaussian o.s. $\cF$ is comprised of affine functions $\phi(\omega)=\phi_0+\sum_{i=1}^n\phi_i\omega_i=:\phi_0+\uphi^T\omega$ on the observation space (which now is $\bR^m$), and, as is immediately seen,
$$
\Phi_\cO(\phi;\mu)=\phi_0+\uphi^T\mu+\half\uphi^T\uphi:(\bR\times\bR^m)\times\bR^m\to\bR.
$$
A straightforward computation shows that in the case in question, using the notation
$\theta=\ln(2I/\epsilon)$,
we get
\begin{equation}\label{example_1}
\begin{array}{rcl}
\Psi_{i,+}(\alpha,\phi)&=&K\alpha\left[\phi_0/\alpha+\uphi^Ty_i/\alpha+{1\over 2}\uphi^T\uphi/\alpha^2\right]+\alpha\theta-g^Tx_i
=K\phi_0+K\uphi^Ty_i-g^Tx_i+{K\over 2\alpha}\uphi^T\uphi+\alpha\theta\\
\Psi_{j,-}(\alpha,\phi)&=&-K\phi_0-K\uphi^Ty_j+g^Tx_j+{K\over 2\alpha}\uphi^T\uphi+\alpha\theta\\
\Opt_{ij}&=&\inf_{\alpha>0,\phi}{1\over 2}\left[\Psi_{i,+}(\alpha,\phi)+\Psi_{j,-}(\alpha,\phi)\right]\\
&=&{1\over 2}g^T[x_j-x_i]+\inf_{\uphi\in\bR^m}\left[{K\over2}\uphi^T[y_i-y_j]+\inf_{\alpha>0}\left[{K\over2\alpha}\uphi^T\uphi+\alpha\theta\right]\right]\\
&=&{1\over 2}g^T[x_j-x_i]+\inf_\uphi\left[{K\over 2}\uphi^T[y_i-y_j]+\sqrt{2K\theta}\|\uphi\|_2\right]\\
&=&\left\{\begin{array}{ll}{1\over 2} g^T[x_j-x_i],&\|y_i-y_j\|_2\leq 2\sqrt{2\theta/K},\\
-\infty,&\|y_i-y_j\|_2> 2\sqrt{2\theta/K}.
\end{array}\right.
\end{array}
\end{equation}
We see that we can safely set $\phi_0=0$, and that
setting $$\I=\{(i,j):\|y_i-y_j\|_2\leq 2\sqrt{2\theta/K}\},$$ $\Opt_{ij}(K)$ is finite when $(i,j)\in\I$ and is $-\infty$ otherwise; in both cases,
the optimization problem specifying $\Opt_{ij}$ has no optimal solution. Indeed, this clearly is the case when $(i,j)\not\in\I$; when $(i,j)\in\I$, a minimizing sequence is, e.g., $\phi\equiv 0,\alpha_i\to0$, but its limit is not in the minimization domain (on this domain, $\alpha$ should be positive). \footnote{Dealing with this case was exactly the reason why in our construction we required from $\phi_{ij},\alpha_{ij}$ to be feasible, and not necessary optimal, solutions to the optimization problems in question.} In the {considered example}, the simplest way to overcome the difficulty is to restrict the optimization domain $\cF^+$ in (\ref{Optij}) with its compact subset $\{\alpha\geq 1/R,\phi_0=0,\|\uphi\|_2\leq R\}$ with a large $R$ (e.g. $R=10^{20}$). {Therefore}, we specify the entities participating in
(\ref{gij}) as
\begin{equation}\label{example_2}
\begin{array}{c}
\begin{array}{rcl}
\phi_{ij}(\omega)&=&\uphi_{ij}^T\omega,\,\,\uphi_{ij}=\left\{\begin{array}{ll}0,&(i,j)\in \I\\
-R[y_i-y_j]/\|y_i-y_j\|_2,&(i,j)\not\in\I\\
\end{array}
\right.\\
\end{array},\,\,
\begin{array}{rcl}
\alpha_{ij}&=&\left\{\begin{array}{ll}1/R,&(i,j)\in \I\\
\sqrt{{K\over2\theta}}R,&(i,j)\not\in\I\\
\end{array}
\right.\\
\end{array}
\end{array}
\end{equation}
resulting in
\begin{equation}\label{example_3}
\begin{array}{rcl}
\varkappa_{ij}&=&{1\over 2}\left[\Psi_{j,-}(\alpha_{ij},\phi_{ij})-\Psi_{i,+}(\alpha_{ij},\phi_{ij})\right]
={1\over 2}g^T[x_i+x_j]-{K\over 2}\uphi_{ij}^T[y_i+y_j]\\
\rho_{ij}&=&{1\over 2}\left[\Psi_{i,+}(\alpha_{ij},\phi_{ij})+\Psi_{j,-}(\alpha_{ij},\phi_{ij})\right]
={K\over 2\alpha_{ij}}\uphi_{ij}^T\uphi_{ij}+\alpha_{ij}\theta+{1\over 2}g^T[x_j-x_i]+{K\over 2}\uphi_{ij}^T[y_i-y_j]\\
&=&\left\{\begin{array}{ll}{1\over 2}g^T[x_j-x_i]+R^{-1}\theta,&(i,j)\in \I\\
{1\over 2}g^T[x_j-x_i]+[\sqrt{2K\theta}-{K\over 2}\|y_i-y_j\|_2]R,&(i,j)\not\in\I\\
\end{array}\right.\\
\end{array}
\end{equation}
{In the numerical experiments we {report below} we use $n=20$, $m=10$, and $I=100$, with $x_i$, $i\leq I$, drawn independently of each other from $\cN(0,I_n)$, and $y_i=Ax_i$ with randomly generated matrix $A$ (namely, matrix with independent $\cN(0,1)$ entries normalized to have unit spectral norm).
The linear form to be recovered is the first coordinate of $x$, the confidence parameter is set to $\epsilon=0.01$, and $R=10^{20}$. Results of a typical experiment are presented in Figure \ref{fig:risks1}.
\begin{figure}[h!]
\includegraphics[scale=0.5]{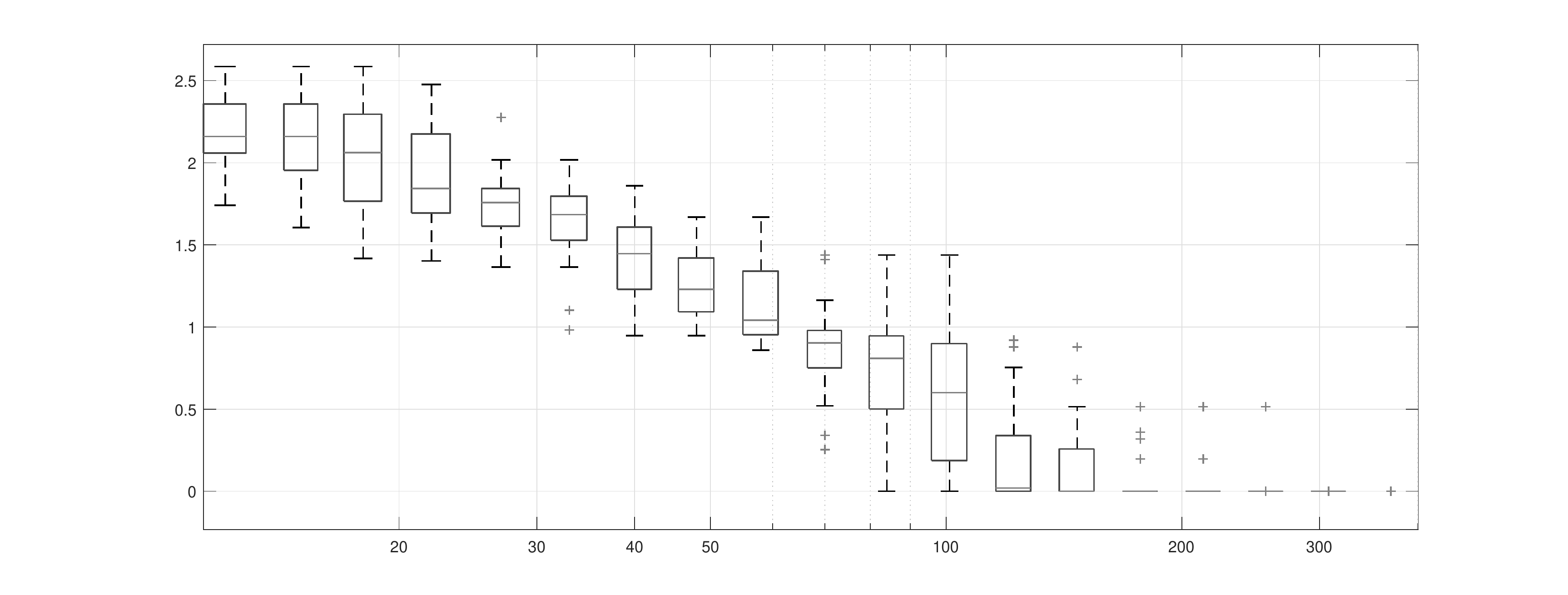}
\caption{\small Boxplot of empirical distributions, over 20 random estimation problems, of the upper 0.01-risk bounds $\max_{1\leq i,j\leq 100}\rho_{ij}$ (as in \rf{example_3}) for different observation sample sizes $K$. }
\label{fig:risks1}
\end{figure}
}

\section{Recovering $N$-convex functions on unions of convex sets}\label{sect:recNconvex}
\subsection{Preliminaries: testing convex hypotheses in good o.s.}\label{sectinfcol}
What follows is a summary of  results of \cite{GJN2015} which are relevant to our current needs.\\
Assume that $\omega^K=(\omega_1,...,\omega_K)$ is a stationary $K$-repeated observation in a good o.s. $\cO=((\Omega,P),\{p_\mu:\mu\in\cM\},\cF)$, so that $\omega_1,...,\omega_K$ are, independently of each other, drawn from a distribution $p_\mu$ with some $\mu\in\cM$. Given $\omega^K$ we want to decide on the hypotheses $H_1$ and $H_2$,
 with $H_\chi$, $\chi=1,2$,  stating that $\omega_t\sim p_\mu$ for some $\mu\in M_\chi$,  where $M_\chi$ is a nonempty convex compact subset of $\cM$.
 In the sequel, we refer to hypotheses of this type, parameterized  by nonempty convex compact subsets of $\cM$,
 as to {\sl convex} hypotheses in the good o.s. in question.\par
The principal ``building block'' of our subsequent constructions is a test $\cT^K$ for this problem which is as follows:
\begin{itemize}
\item Given convex compact sets $M_\chi$, $\chi=1,2$, we solve the optimization problem
\begin{equation}\label{eq1}
\Opt=\max_{\mu\in M_1,\nu\in M_2}\ln\left(\int_\Omega\sqrt{p_\mu(\omega)p_\nu(\omega)}P(d\omega)\right)
\end{equation}
It is shown in \cite{GJN2015} that in the case of good o.s., problem (\ref{eq1}) is a convex problem (convexity meaning that the objective to be maximized is a concave continuous function of $\mu,\nu$) and an optimal solution exists.
\begin{quote}
\noindent Note that for basic good o.s.'s problem (\ref{eq1}) reads
\begin{equation}\label{eq11}
\Opt=\max_{\mu\in M_1,\nu\in M_2}\left\{\begin{array}{ll}
-\hbox{\small $1\over 8$}\|\mu-\nu\|_2^2,&\hbox{Gaussian o.s.}\\
-{\half}\sum_{i=1}^d [\sqrt{\mu_i}-\sqrt{\nu_i}]^2,&\hbox{Poisson o.s.}\\
\ln\left(\sum_{i=1}^d\sqrt{\mu_i\nu_i}\right),&\hbox{Discrete o.s.}\\
\end{array}\right.
\end{equation}
\end{quote}
\item An optimal solution $\mu_*$, $\nu_*$ to (\ref{eq1}) induces  {\sl detectors}
\begin{equation}\label{eq2}
\begin{array}{rcl}\phi_*(\omega)&=&{\half}\ln(p_{\mu_*}(\omega)/p_{\nu_*}(\omega)):\;\Omega\to\bR,
\\\phi_*^{(K)}(\omega^K)&=&\sum_{t=1}^K\phi_*(\omega_t):\;\; \Omega\times...\times\Omega\to\bR
\end{array}\end{equation}
Given a stationary $K$-repeated observation $\omega^K$, the test $\cT^K$ accepts hypothesis $H_1$ and rejects hypothesis $H_2$ whenever $\phi_*^{(K)}(\omega^K) \geq0$, otherwise the test rejects $H_1$ and accepts $H_2$. The {\sl risk} of $\cT^K$ -- the maximal probability to reject a hypothesis when it is true -- does not exceed $\epsilon_\star^K,$ where
$$\epsilon_\star=\exp(\Opt).
$$
In other words, whenever observation $\omega^K$ stems from a distribution $p_\mu$ with $\mu\in M_1\cup M_2$,
\begin{itemize}
\item the $p_\mu$-probability to reject $H_1$ when the hypothesis is true (i.e., when $\mu\in M_1$) is at most $\epsilon_\star^K$,     and
\item the $p_\mu$-probability to reject $H_2$ when the hypothesis is true (i.e., when $\mu\in M_2$) is at most $\epsilon_\star^K$.
\end{itemize}
\end{itemize}
The test $\cT^K$ possesses the following optimality properties:
\begin{itemize}
\item[{\bf A.}] The associated detector {\sl $\phi_*^{(K)} $ and the risk $\epsilon_\star^K$ form an optimal solution and the optimal value in the optimization problem
\[
\begin{array}{c}
\min\limits_\phi\max\left[\max_{\mu\in M_1}\int_{\Omega^K} \e^{-\phi(\omega^K)}p^{(K)}_\mu(\omega^K)P^K(d\omega^K),
\max_{\nu\in M_2}\int_{\Omega^K} \e^{\phi(\omega^K)}p^{(K)}_\nu(\omega^K)P^K(d\omega^K)\right],\\
\big[\Omega^K=\underbrace{\Omega\times...\times\Omega}_{K},\;\;p_\mu^{(K)}(\omega^K)=\prod_{t=1}^Kp_\mu(\omega_t),\,
P^K=\underbrace{P\times...\times P}_{K}\big]\\
\end{array}
\]
where the minimum is taken w.r.t. all Borel functions $\phi(\cdot):\Omega^K\to\bR$;}
\item[{\bf B.}] {\sl Let $\epsilon\in (0,1/2)$, and suppose that there exists a test which, using a stationary $\overline{K}$-repeated observation, decides on the hypotheses $H_1$, $H_2$ with risk $\leq\epsilon$. Then
    \begin{equation}\label{epsilonstar}
    \epsilon_\star\leq [2\sqrt{\epsilon(1-\epsilon)}]^{1/\overline{K}}
    \end{equation}
    and the test $\cT^K$ with
    $$
    K=\left\rfloor {2\ln(1/\epsilon)\over \ln\left( [4\epsilon(1-\epsilon)]^{-1}\right)} \overline{K}\right\lfloor
    $$
    decides on the hypotheses $H_1,H_2$ with risk $\leq\epsilon$ as well. Note that $K=2(1+o(1))\overline{K}$ as $\epsilon\to+0$.}
\end{itemize}
\paragraph{``Inferring colors:'' testing multiple hypotheses in good o.s.} As shown in \cite{GJN2015}, the just outlined near-optimal pairwise tests deciding on pairs of convex hypotheses in good o.s.'s can be used as building blocks when constructing near-optimal tests deciding on multiple convex hypotheses. In the sequel, we will repeatedly use one of these constructions, namely, as follows.
\par
Assume that we are given a good o.s. $\cO=((\Omega,P),\{p_\mu:\mu\in\cM\},\cF)$ and two finite collections of nonempty convex compact subsets  $B_1,...,B_b$ (``blue sets'') and $R_1,...,R_r$ (``red sets'')   of $\cM$. Our objective is, given a stationary $K$-repeated observation $\omega^K$ stemming from a distribution $p_\mu$, $\mu\in\cM$, to infer the color of $\mu$, that is, to decide on the hypothesis $\mu\in B:=B_1\cup...\cup B_b$    vs. the alternative $\mu\in R:=R_1\cup...\cup R_r$. To this end we act as follows:
\begin{enumerate}
\item For every pair $i,j$ with $i\leq b$ and $j\leq r$, we solve the problem (\ref{eq11})  with $B_i$ in the role of $M_1$ and $R_j$ in the role of $M_2$; we denote $\Opt_{ij}$ the associated optimal values. The corresponding optimal solutions $\mu_{ij}$ and $\nu_{ij}$ give rise to the detectors
\begin{equation}\label{eqij}
\begin{array}{c}
\phi_{ij}(\omega)={1\over 2}\ln\left(p_{\mu_{ij}}(\omega)/p_{\nu_{ij}}(\omega)\right):\Omega\to\bR,\,\,
\phi_{ij}^{(K)}=\sum_{t=1}^K\phi_{ij}(\omega_t):\Omega^K\to\bR\\
\end{array}
\end{equation}
(cf. (\ref{eq2})) and {\sl risks}
\begin{equation}\label{epsij}
\epsilon_{ij}=\exp(\Opt_{ij})=\int_\Omega \sqrt{p_{\mu_{ij}}(\omega)p_{\nu_{ij}}(\omega)}P(d\omega).
\end{equation}
\item We build the entrywise positive $b\times r$ matrix $E^{(K)}=[\epsilon_{ij}^K]_{{1\leq i\leq b\atop 1\leq j\leq r}}$ and symmetric entrywise nonnegative $(b+r)\times (b+r)$ matrix $E_K=\hbox{\tiny$\left[\begin{array}{c|c}&E^{(K)}\cr\hline [E^{(K)}]^T&\cr\end{array}\right]$}$. Let $\epsilon_{K}$ be the spectral norm of the matrix $E^{(K)}$ (equivalently, spectral norm of $E_K$), and let $e=[g;h]$\footnote{We use ``Matlab notation'' $[a;b]$ for vertical and $[a,b]$ for horizontal concatenation of matrices $a, b$ of appropriate dimensions.} be the Perron-Frobenius eigenvector of $E_K$, so that $e$ is a nontrivial nonnegative vector such that $E_Ke=\epsilon_Ke$. Note that from entrywise positivity of $E^{(K)}$ it immediately follows that $e>0$, so that the quantities
    $$
    \alpha_{ij}=\ln(h_j/g_i),\;\;1\leq i\leq b,\; 1\leq j\leq r
    $$
    are well defined. We set
    \begin{equation}\label{eq12}
    \psi_{ij}^{(K)}(\omega^K)=\phi_{ij}^{(K)}(\omega^K)-\alpha_{ij}=\sum_{t=1}^K\phi_{ij}(\omega_t)-\alpha_{ij}:\Omega^K\to\bR,\,1\leq i\leq b,1\leq j\leq r
\end{equation}
\item Given observation $\omega^K\in\Omega^K$ with $\omega_t$, $t=1,...,K$, drawn, independently of each other, from a distribution $p_\mu$, we
claim that $\mu$ is blue (equivalently, $\mu\in B$), if there exists $i\leq b$ such that $\psi_{ij}(\omega^K)\geq0$ for all $j=1,...,r$, and claim that $\mu$ is red (equivalently, $\mu\in R$) otherwise.
\end{enumerate}
The main result about the just described ``color inferring'' test is as follows
\begin{proposition}\label{basicprop} {\rm \cite[Proposition 3.2]{GJN2015}}  Let the components  $\omega_t$ of $\omega^K$ be drawn, independently of each other, from distribution $p_\mu\in B\cup R$. Then the just defined test, for every $\omega^K$, assigns $\mu$ with exactly one color, blue or red, depending on the observation. Moreover,
\begin{itemize}
\item when $\mu$ is blue (i.e., $\mu\in B$), the test makes correct inference ``$\mu$ is blue'' with $p_\mu$-probability at least $1-\epsilon_K$;
\item similarly, when $\mu$ is red (i.e., $\mu\in \bR$), the test makes correct inference ``$\mu$ is red'' with $p_\mu$-probability at least $1-\epsilon_K$.
\end{itemize}
\end{proposition}

\subsection{Problem's setting}\label{sect:setting}
In the sequel, we deal with the situation as follows. Given are:
\begin{enumerate}
\item good o.s.  $\cO=((\Omega,P),\{p_\mu(\cdot):\mu\in\cM\},\cF)$,
\item convex compact set $\cX\subset\bR^n$ along with a collection of $I$ convex compact sets $X_i\subset\cX$,
\item affine ``encoding'' $x\mapsto A(x): \cX\to\cM$,
\item a continuous function $f(x):\cX\to\bR$ which is {\sl $N$-convex}, meaning that for every $a\in\bR$ the sets
$\cX^{a,\geq}=\{x\in\cX: f(x)\geq a\}$ and $\cX^{a,\leq}=\{x\in\cX: f(x)\leq a\}$ can be represented as unions of at most $N$ closed convex sets $\cX^{a,\geq}_{\nu}$, $\cX^{a,\leq}_{\nu}$:
\begin{equation}\label{Meq1}
\cX^{a,\geq}=\bigcup\limits_{\nu=1}^N\cX^{a,\geq}_{\nu},\,\,\cX^{a,\leq}=\bigcup\limits_{\nu=1}^N\cX^{a,\leq}_{\nu}.
\end{equation}
\end{enumerate}
For some {\sl unknown} $x$ known to belong to $X=\bigcup\limits_{i=1}^IX_i$, we have at our disposal observation $\omega^K=(\omega_1,...,\omega_K)$ with i.i.d. $\omega_t\sim p_{A(x)}(\cdot)$, and our goal is to estimate from this observation the quantity $f(x)$.
\par
The {\em $\epsilon$-risk} of a candidate estimate $\widehat{f}(\omega^K)$ is defined in the same way it was done in Section \ref{sect:linear:prob}. Specifically, given tolerances  $\rho>0$, $\epsilon\in(0,1)$, we call $\widehat{f}(\omega^K)$ {\sl $(\rho,\epsilon)$-reliable}, if for every $x\in X$, $|\widehat{f}(\omega^K)-f(x)|\leq\rho$ with the $p_{A(x)}$-probability at least $1-\epsilon$.
The $\epsilon$-risk of $\widehat{f}(\omega^K)$ is the smallest $\rho$ such that  $\widehat{f}(\cdot)$ is $(\rho,\epsilon)$-reliable.
%
%
\paragraph{Examples of $N$-convex functions.}
In the above problem setting we allow 
$X$ to be a finite union of convex sets, and function $f$ is assumed to be $N$-convex. Being rather restrictive, the latter class comprises, along with linear functions, some interesting examples, which we discuss below.
\begin{example}\label{ExMinMax}{\rm[Minima and Maxima of linear-fractional functions]} Every function which can be obtained from linear-fractional functions ${g_{\nu}(x)\over h_{\nu}(x)}$ ($g_{\nu}$, $h_{\nu}$ are affine functions on $\cX$, and $h_{\nu}$ are positive on $\cX$) by taking maxima and minima is $N$-convex for appropriately selected $N$due to the following immediate observations:
\begin{itemize}
\item linear-fractional function ${g(x)\over h(x)}$ with a denominator which is positive on $\cX$ is 1-convex;
\item if $f(x)$ is $N$-convex, so is $-f(x)$;
\item if $f_i(x)$ is $N_i$-convex, $i=1,2,...,I$, then $f(x)=\max_if_i(x)$ is $\max[\prod_iN_i,\sum_iN_i]$-convex.
\par
{\rm Indeed, we have
$$
\{x\in\cX:f(x)\leq a\}=\bigcap\limits_{i=1}^I\{x\in\cX:f_i(x)\leq a\},\;\mbox{and}\;\{x\in\cX:f(x)\geq a\}=\bigcup\limits_{i=1}^I\{x\in\cX:f_i(x)\geq a\}.
$$
The first set is the intersection of $I$ unions of convex sets with $N_i$ components in $i$-th union, and thus is the union of $\prod_iN_i$ convex sets; the second  set
is the union of $I$ unions, $N_i$ components in the $i$-th of them, of convex sets, and thus is the union of $\sum_iN_i$ convex sets.}
\end{itemize}
\end{example}
\begin{example}\label{ExCondQuant} {\em [Conditional quantile] Let $S=\{s_1<s_2<...<s_M\}\subset \bR$. For a nonvanishing probability distribution $q$ on $S$ and $\alpha\in[0,1]$, let
$\chi_\alpha[q]$ be the {\sl regularized} $\alpha$-quantile of $q$ defined as follows: we pass from $q$ to the distribution on $[s_1,s_M]$ by spreading uniformly the mass $q_{\nu}$, $\aic{1}{2}<\nu\leq M$, over $[s_{\nu-1},s_{\nu}]$, and assigning  mass $q_1$ to the point $s_1$;
$\chi_\alpha[q]$ is the usual $\alpha$-quantile
of the resulting distribution $\bar{q}$:
{
$$
\begin{array}{c}
\chi_\alpha[q]=\min\{s\in[s_1,s_M]: \bar{q}\{[s_1,s]\}\geq\alpha\}.\\
\epsfxsize=200pt\epsfysize=150pt\epsffile{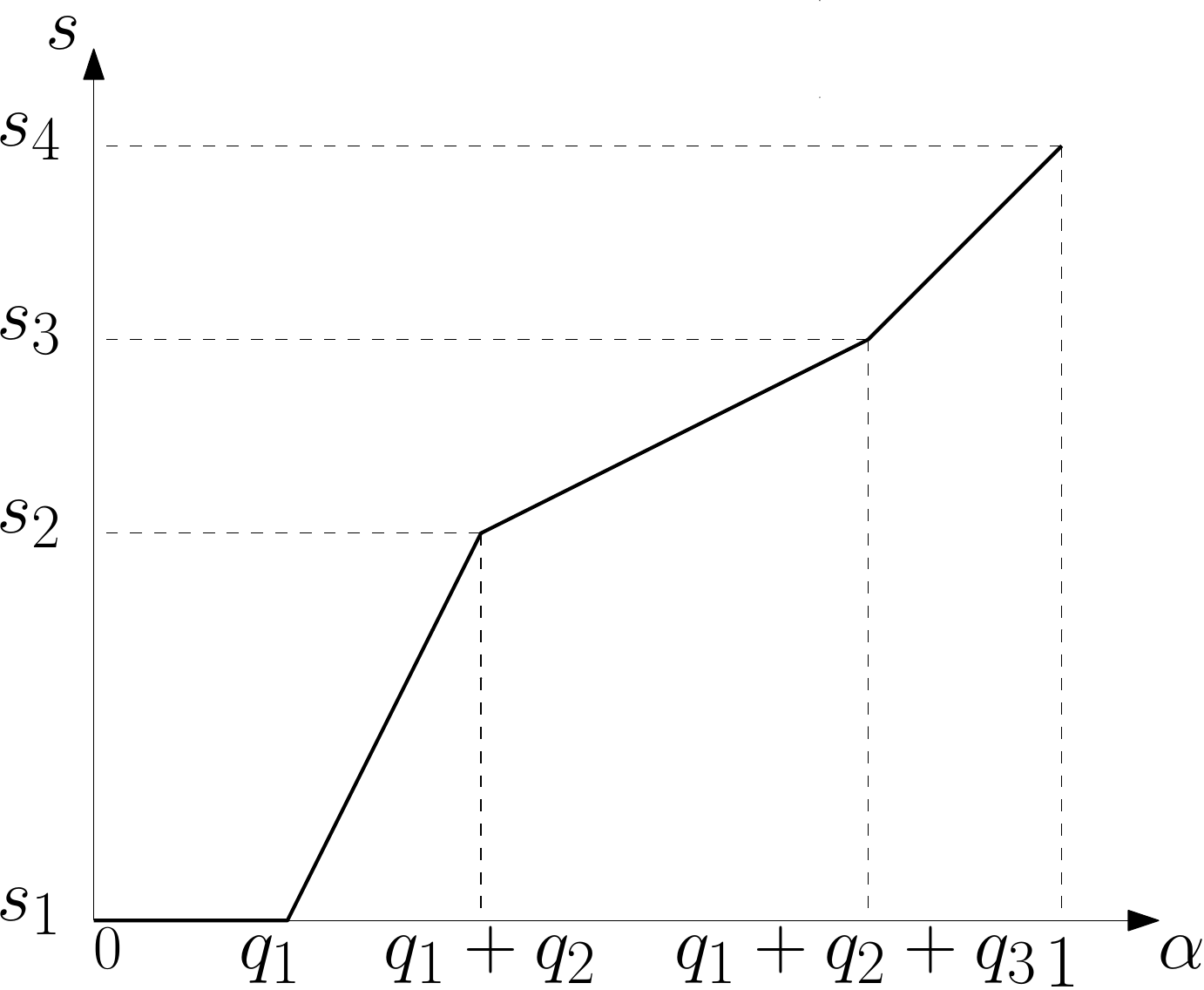}\\
\hbox{Regularized quantile as function of $\alpha$, $M=4$}\\
\end{array}
$$
}
\vskip5pt
Given, along with $S$, a finite set $T$, let $\cX$ be a convex compact set in the space of nonvanishing probability distributions on $S\times T$.
Given $\tau\in T$, consider the conditional, by the condition $t=\tau$, distribution $p_\tau(\cdot)$ of $s\in S$ induced by a distribution $p(\cdot,\cdot)\in \cX$:
$$
p_\tau(\mu)={p(\mu,\tau)\over\sum_{\nu=1}^Mp(\nu,\tau)},\,1\leq\mu\leq M,
$$
where $p(\mu,\tau)$ is the $p$-probability for $(s,t)$ to take value $(s_\mu,\tau)$, and $p_\tau(\mu)$ is the $p_\tau$-probability for $s$ to take value $s_\mu$, $1\leq\mu\leq M$. \par
The function $\chi_\alpha[p_\tau]:\cX\to\bR$ turns out to be 1-convex, see Appendix \ref{1convex}.}
\end{example}

\subsection{Bisection Estimate}
As we have already mentioned, the proposed estimation procedure is a ``close relative'' of the binary search algorithm of \cite{donoho1991geom2}, but is not identical to that algorithm.
Though the bisection estimator is, in a nutshell, quite simple, its formal description turns out to be rather involved. For this reason we start its presentation with an informal outline, which exposes some simple ideas underlying the construction.
\subsubsection{Outline}\label{outline}

\begin{figure}
\centering
\begin{tabular}{ccc}
\includegraphics[scale=0.3]{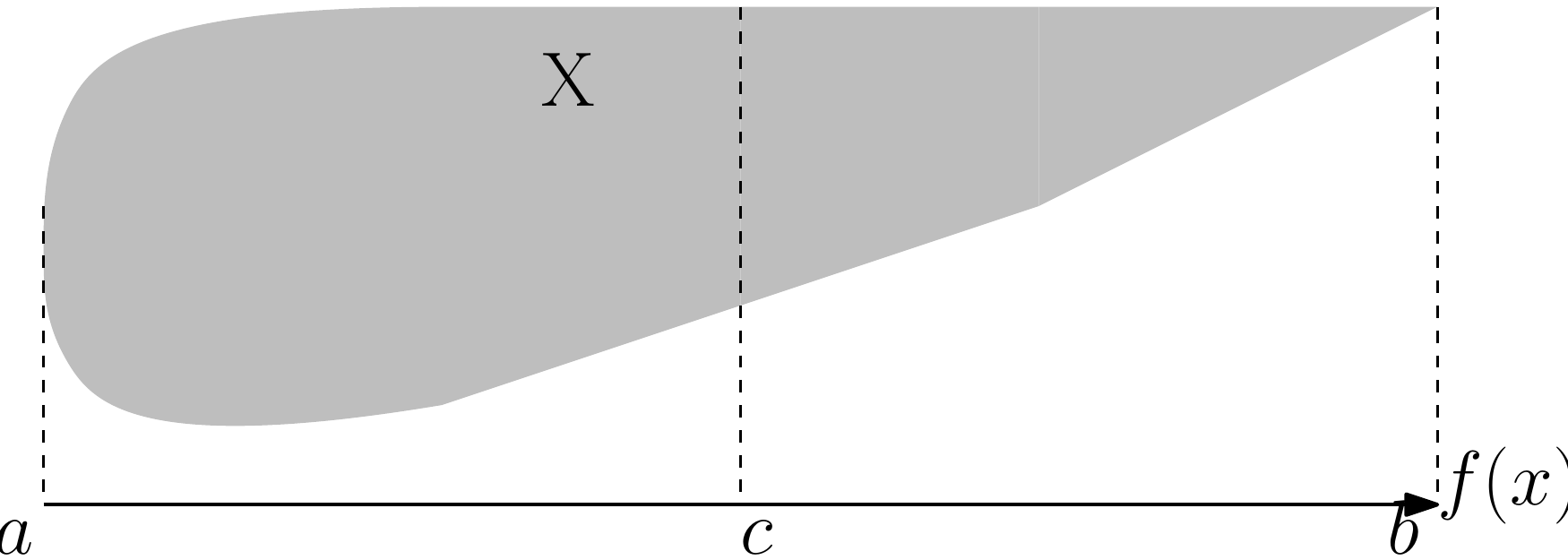}
&
\includegraphics[scale=0.3]{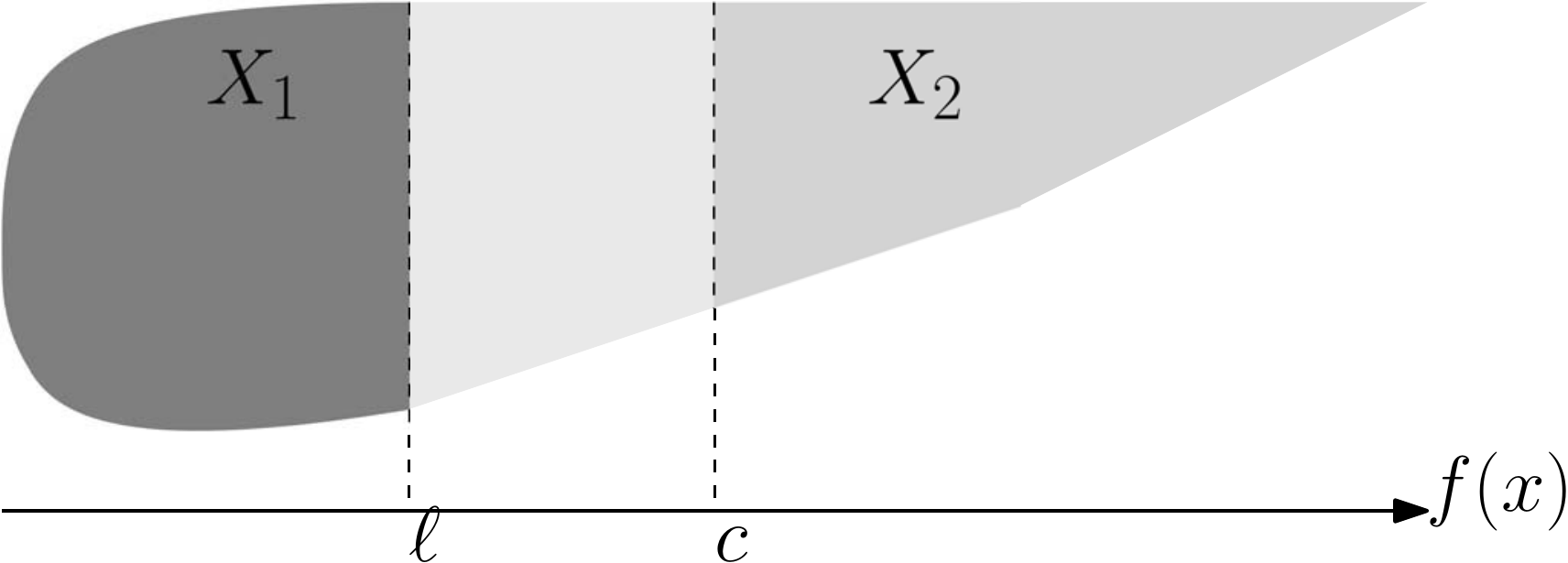}
&
\includegraphics[scale=0.3]{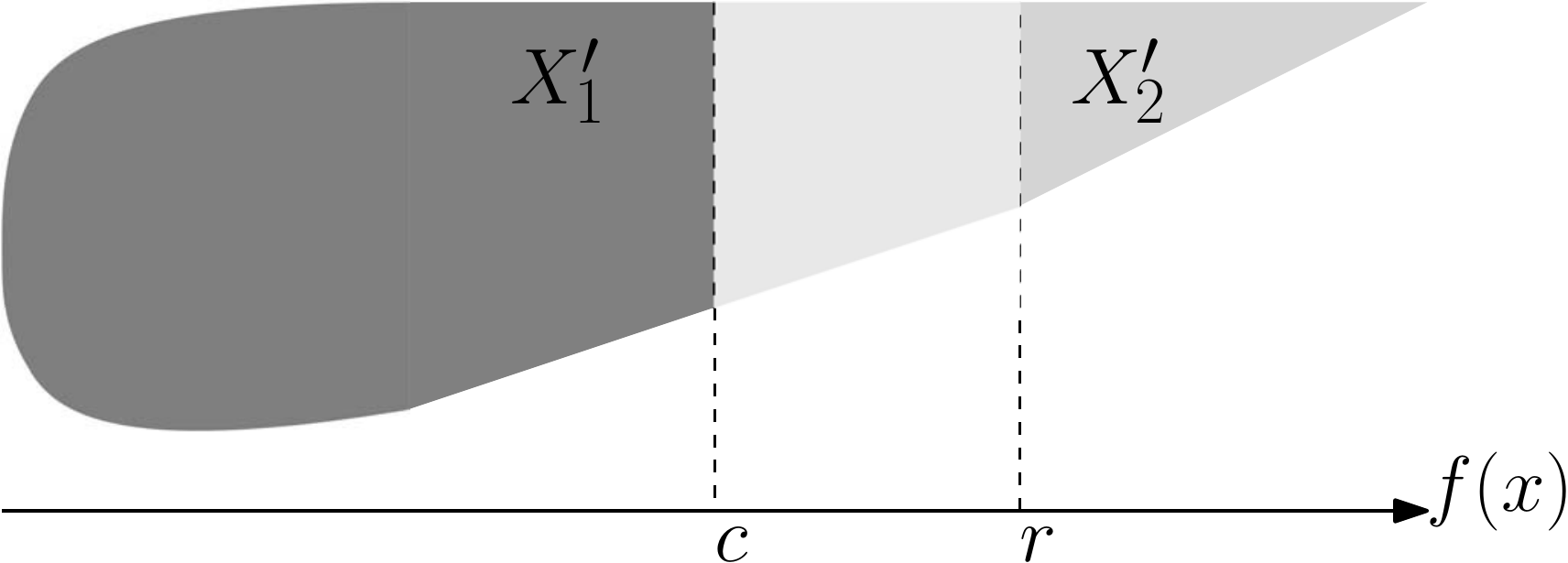}\\
$(a)$&$(b)$&$(c)$\\~
\end{tabular}
\caption{\label{algopic}  {\small Bisection via hypothesis testing. $(a)$: set $X$ of signals and initial localizer $[a,b]$ for the value of $f(x)=x_1$;
$(b)$: left hypothesis $H_1=\{x\in X_1\}$ and
right hypothesis $H_2=\{x\in X_2\}$;
$(c)$: left hypotheses $H_1'=\{x\in X_1'\}$ and right hypothesis $H_2'=\{x\in X_2'\}$.
}}
\end{figure}
\par
Let us consider a simple situation where the signal space $X$ is a convex set in $\bR^2$, as presented in Figure \ref{algopic}, and suppose that our objective is to estimate the value of a linear function $f(x)=x_1$ at $x=[x_1;x_2]\in X$ given a Gaussian observation $\omega$ with mean $A(x)$, where $A(\cdot)$ is a given affine mapping, and known covariance. Observe that hypotheses $f(x)\geq b$ and $f(x)\leq a$ translate into convex hypotheses on the expectation of the observed Gaussian r.v., so that we can use the hypothesis testing machinery of Section \ref{sectinfcol} to decide on hypotheses of this type and to localize $f(x)$ in a (hopefully, small) segment by a bisection-type process. Before describing the process, let us make a terminological agreement. In the sequel we sometimes use pairwise hypothesis tests in the situation where {\em neither} of the hypotheses is true. In this case, we say that the outcome of a test is correct, if the rejected hypothesis indeed is wrong; in this case, the accepted hypothesis can be wrong as well, but this can happen only when both tested hypotheses are wrong.
\par
Let $\epsilon \in (0,1)$ and let $L$ be a positive integer. The estimation procedure is organized in steps. At the beginning of the first step $\Delta_1=[a,b]$ with $a=\min_{x\in X} x_1$, and  $b=\max_{x\in X} x_1$, is the current localizer for the value of $f(x)=x_1$, see Figure \ref{algopic}, and let $c=\half(a+b)$.
To compute the new localizer, we run {\em a pair of Left vs. Right tests} $\cT$ and $\cT'$, such that
\begin{itemize}
\item $\cT$ decides upon the ``left pair of hypotheses'' $H_1=\{x\in X:\,x_1\leq \ell\}$ (left) vs. $H_2=\{x\in X:\,x_1\geq c\}$ (right), where $\ell<c$ is as close to $c$ as possible under the restriction that $\cT$ decides on $H_1$, $H_2$ with risk $\leq {\epsilon\over 2L}$;
\item $\cT'$ decides upon the ``right pair of hypotheses'' $H^\prime_1=\{x\in X:\,x_1\leq c\}$ (left) vs. $H^\prime_2=\{x\in X:\,x_1\geq r\}$ (right), where $r>c$ is as close to $c$ as possible under the restriction that $\cT'$ decides on $H_1^\prime,H_2^\prime$ with risk $\leq {\epsilon\over 2L}$.
\end{itemize}
Assuming that both tests rejected wrong hypotheses (this happens with probability at least $1-{\epsilon\over L}$), the results of the tests allow for the following conclusions:
\begin{itemize}
\item  when both tests reject right hypotheses from the corresponding pairs, it is certain that $x_1\leq c$  (since otherwise in the first test the rejected hypothesis were in fact true, contradicting the assumption that both tests make no wrong rejections);
\item when both tests reject left hypotheses from the corresponding pairs, it is certain that  $x_1\geq c$ (for the same reasons as in the previous case);
\item when the tests ``disagree,'' rejecting hypotheses of different colors, $x_1\in [\ell,r]$. Indeed, otherwise either $x_1\leq \ell$ (and thus $x$ is ``colored left'' in both pairs of hypotheses), or $x_1\geq r$ (and $x$ is ``colored right'' in both pairs). Since we have assumed that in both tests no wrong rejections took place, in the first case both tests must reject right hypotheses, and both should reject left ones in the second, while none of these events took place.
\end{itemize}
In the first two cases we take the right or the left half of the initial segment $\Delta_1=[a,b]$ as a new localizer for $f(x)=x_1$ (and the corresponding cut $X\cap \{x_1\geq c\}$ or $X\cap \{x_1\leq c\}$ as a new localizer for $x$). In the last case, we take the segment $[\ell,r]$ as a new localizer for $x_1$, {\em terminate the process and output $\wh{f}=\half(\ell+r)$ as estimate of $f(x)$} -- the $\epsilon/ L$-risk of this estimate is equal to $\half(r-\ell)$ and is already small!
In Bisection, we iterate the outlined procedure, replacing current localizers with twice smaller ones until terminating either due to running into ``disagreement,''  or due to reaching a prescribed number $L$ of steps. Upon termination, we return the last localizer as a confidence set for
$f(x)=x_1$, and its midpoint -- as the estimate of $f(x)$.\par
Note that, unlike the binary search procedure of \cite{donoho1991geom2}, in our procedure the ``search trajectory'' -- the sequence of pairs of hypotheses participating in the tests -- is not random, it is uniquely defined by the value of $f(x)$, provided no wrong rejections happen. Indeed, with no wrong rejections prior to termination, {the sequence of localizers produced by the procedure is {\em exactly {the same} as if} we were running deterministic bisection algorithm, that is, were updating  subsequent localizers $\Delta_\ell$ for $f(x)$ according to the rules
\begin{itemize}
\item $\Delta_1=[a,b]$, the obvious initial segment $f(x)$,
\item $\Delta_{\ell+1}$ is precisely the half of $\Delta_\ell$ containing $f(x)$ (say, the left half in the case of a tie).
\end{itemize}
In the above argument we neglected the possibility of wrong rejection by one of the tests we ran. Since, by construction, the risk of each test does not exceed ${\epsilon\over 2L}$ and, by the above,  with no wrong rejections, {\sl the sequence of tests we run depends solely on the value $f(x)$, not on the observations} (observations can affect only the number of steps before termination), the probability of wrong rejection in course of running the algorithm is
$\leq \epsilon$.
Note that the  risks of ``individual tests'' define, in turn, the allowed width of separators -- segments $[\ell,c]$ and $[c,r]$ in Figure \ref{algopic}.b (``uncertainty zone'' of the corresponding test), and thus -- the accuracy to which $f(x)$ can be estimated. It should be noted that the number $L$ of steps of Bisection always is a moderate integer. Indeed, otherwise the width of the separators at the concluding bisection steps (which is of order of $2^{-L}$), would be too small to allow for deciding on the concluding pairs of our hypotheses with risk ${\epsilon\over 2L}$.
\par
From the above sketch of our construction, it is clear that all that matters is our ability to decide, given $\ell<r$,
on the pairs of hypotheses $\{x\in X: \,f(x)\leq \ell\}$ and $\{x\in X: \,f(x)\geq r\}$ via observation drawn from $p_{A(x)}$. In our outline, these were convex hypotheses in Gaussian o.s., and in this case we can use detector-based pairwise tests presented in Section \ref{sectinfcol}. Applying the machinery developed in the latter section, we could also handle the case when the sets $\{x\in X: f(x)\leq \ell\}$ and $\{x\in X: f(X)\geq r\}$ are finite unions of convex sets (which is the case when
$f$ is $N$-convex and $X$ is a finite union of convex sets), the o.s. in question still being good, and this is the situation we intend to consider.
\subsubsection{Building the Bisection estimate: preliminaries}
While the construction we present below admits numerous refinements, we focus here on its simplest version as follows (for notation, see Section  \ref{sect:setting}).

\paragraph{Upper and lower feasibility/infeasibility, sets $Z^{a,\geq}_i$ and $Z^{a,\leq}_i$.}
Let $a$ be a real. We associate with $a$ the collection of  {\sl upper $a$-sets} defined as follows: we look at the sets $X_i\cap \cX^{a,\geq}_{\nu}$, $1\leq i\leq I$, $1\leq \nu\leq N$, and arrange the nonempty sets from this family into a sequence $Z^{a,\geq}_i$, $1\leq i\leq I_{a,\geq}$. Here $I_{a,\geq}=0$ if all sets in the family are empty; in the latter case, we refer to $a$ as {\sl upper-infeasible}, and {\sl upper-feasible} otherwise. Similarly, we associate with $a$ the collection of {\sl lower $a$-sets} $Z^{a,\leq}_i$, $1\leq i\leq I_{a,\leq}$ by arranging into a sequence all nonempty sets from the family $X_i\cap \cX^{a,\leq}_\nu$, $1\leq i\leq I$, $1\leq \nu\leq N$. We say that $a$ is lower-feasible or lower-infeasible depending on whether
$I_{a,\leq}$ is positive or zero. Note that upper and lower $a$-sets, if any, are nonempty convex compact sets, and
\begin{equation}\label{upperlowersets}
X^{a,\geq}:=\{x\in X:f(x)\geq a\}=\bigcup\limits_{1\leq i\leq I_{a,\geq}}Z^{a,\geq}_i,\,\,X^{a,\leq}:=\{x\in X:f(x)\leq a\}=\bigcup\limits_{1\leq i\leq I_{a,\leq}}Z^{a,\leq}_i.
\end{equation}
\paragraph{Right tests.} Given a segment $\Delta=[a,b]$ of positive length with lower-feasible $a$, we associate with this segment {\sl right test} -- a function $\cT^K_{\Delta\r}(\omega^K)$ taking values $\up$ and $\dw$, and risk $\sigma_{\Delta\r}\geq0$ -- as follows:
\begin{enumerate}
\item if $b$ is upper-infeasible, $\cT^K_{\Delta\r}(\cdot)\equiv\dw$ and $\sigma_{\Delta\r}=0$.
\item if $b$ is upper-feasible, the collections $\{A(Z^{b,\geq}_i)\}_{i\leq I_{b,\geq}}$ (``right sets''), $\{A(Z^{a,\leq}_j)\}_{j\leq I_{a,\leq}}$ (``left sets''),  are nonempty, and the test is the associated with these sets Inferring Color test from Section \ref{sectinfcol} {\sl as applied to the stationary $K$-repeated version of $\cO$}  {\sl in the role of $\cO$},  specifically,
 \begin{itemize}
 \item for $1\leq i\leq I_{b,\geq}$, $1\leq j\leq I_{a,\leq}$, we build the detectors $\phi_{ij\Delta}^K(\omega^K)=\sum_{t=1}^K \phi_{ij\Delta}(\omega_t)$, with $\phi_{ij\Delta}(\omega)$ given by
 \begin{equation}\label{Neq4r}
 \begin{array}{rcl}
 (r_{ij\Delta},s_{ij\Delta})&\in&\Argmin_{r\in Z^{b,\geq}_i,s\in Z^{a,\leq}_j}\ln\left(\int_\Omega\sqrt{p_{A(r)}(\omega)p_{A(s)}(\omega)}P(d\omega)\right),\\
 \phi_{ij\Delta}(\omega)&=&{1\over 2}\ln\left(p_{A(r_{ij\Delta})}(\omega)/p_{A(s_{ij\Delta})}(\omega)\right)\\
 \end{array}
 \end{equation}
 set
 \begin{equation}\label{Neq5r}
 \epsilon_{ij\Delta}=\int_\Omega\sqrt{p_{A(r_{ij\Delta})}(\omega)p_{A(s_{ij\Delta})}(\omega)}P(d\omega)
 \end{equation}
 and build the $I_{b,\geq}\times I_{a,\leq}$ matrix $E_{\Delta\r}=[\epsilon_{ij\Delta}^K]_{{1\leq i\leq I_{b,\geq}\atop 1\leq j\leq I_{a,\leq}}}$;
 \item $\sigma_{\Delta\r}$ is defined as the spectral norm of $E_{\Delta\r}$. We compute the Perron-Frobenius eigenvector $[g^{\Delta\r};h^{\Delta\r}]$ of the matrix  $\left[\begin{array}{c|c}&E_{\Delta\r}\cr\hline E_{\Delta\r}^T\cr\end{array}\right]$, so we have (see Section \ref{sectinfcol})
 $$
 g^{\Delta\r}>0,\,h^{\Delta\r}>0, \sigma_{\Delta\r} g^{\Delta\r}=E_{\Delta\r}h^{\Delta\r},\,\sigma_{\Delta\r}  h^{\Delta\r}=E_{\Delta\r}^Tg^{\Delta\r}.
 $$
 Finally, we define the matrix-valued function
 $$
 D_{\Delta\r}(\omega^K)=[\phi_{ij\Delta}^K(\omega^K)+\ln(h^{\Delta\r}_j)-\ln(g^{\Delta\r}_i)]_{{1\leq i\leq I_{b,\geq}\atop1\leq j\leq I_{a,\leq}}}.
 $$
 Test $\cT^K_{\Delta\r}(\omega^K)$ takes value $\up$ iff the matrix $D_{\Delta\r}(\omega^K)$ has a nonnegative row, and takes value $\dw$ otherwise.
 \end{itemize}
 \end{enumerate}
Given $\delta>0$ and $\varkappa>0$, we call segment $\Delta=[a,b]$ {\sl $\delta$-good (right)}, if $a$ is lower-feasible, $b>a$, and $\sigma_{\Delta\r}\leq\delta$ and call  a $\delta$-good (right) segment $\Delta=[a,b]$ {\sl $\varkappa$-maximal}, if the segment $[a,b-\varkappa]$ is not $\delta$-good (right).
\paragraph{Left tests.}
The ``mirror'' version of the above is as follows. Given a segment $\Delta=[a,b]$ of positive length with upper-feasible $b$, we associate with this segment {\sl left test} -- a function $\cT^K_{\Delta\l}(\omega^K)$ taking values $\up$ and $\dw$, and risk $\sigma_{\Delta\l}\geq0$ -- as follows:
\begin{enumerate}
\item if $a$ is lower-infeasible, $\cT^K_{\Delta\l}(\cdot)\equiv\up$ and $\sigma_{\Delta\l}=0$.
\item if $a$ is lower-feasible, we set $\cT^K_{\Delta\l}\equiv \cT^K_{\Delta\r}$, $\sigma_{\Delta\l}=\sigma_{\Delta\r}$.
\end{enumerate}
 Given $\delta>0$, $\varkappa>0$, we call segment $\Delta=[a,b]$ {\sl $\delta$-good (left)}, if $b$ is upper-feasible, $b>a$, and $\sigma_{\Delta\l}\leq\delta$ and  call a $\delta$-good (left) segment $\Delta=[a,b]$ {\sl $\varkappa$-maximal}, if the segment $[a+\varkappa,b]$ is not $\delta$-good (left).
\paragraph{Remark:}  note that when $a<b$ and $a$ is lower-feasible, $b$ is upper-feasible, so that the sets
$$
X^{a,\leq}=\{x\in X: f(x)\leq a\},\,X^{b,\geq}=\{x\in X: f(x){\geq b}\}
$$
are nonempty,
the right and the left  tests  $\cT^K_{\Delta\l}$, $\cT^K_{\Delta\r}$ are identical and coincide with the Color Inferring test, built as explained in Section \ref{sectinfcol}, deciding, via stationary $K$-repeated observations,  on the ``type'' of the distribution $p_{A(x)}$ underlying observations -- whether this type is left (``left'' hypothesis stating that $x\in X$ and $f(x)\leq a$,  whence $A(x)\in\bigcup\limits_{1\leq i\leq I_{a,\leq}}A(Z^{a,\leq}_i)$), or right (``right'' hypothesis, stating that $x\in X$ and $f(x)\geq b$, whence  $A(x)\in\bigcup\limits_{1\leq i\leq I_{b,\geq}}A(Z^{b,\geq}_i)$).
When $a$ is lower-feasible and $b$ is {\sl not} upper-feasible, the right hypothesis is empty, and the left test associated with $[a,b]$, naturally, always accepts the left hypothesis. Similarly, when $a$ is lower-infeasible and $b$ is upper-feasible, the right test associated with $[a,b]$ always accepts the right hypothesis. \par
A segment $[a,b]$ with $a<b$ is $\delta$-good (left), if the corresponding to the segment ``right'' hypothesis is nonempty, and the left  test $\cT^K_{\Delta \l}$ associated with $[a,b]$ decides on the ``right'' and the ``left'' hypotheses with risk $\leq\delta$, that is,
\begin{itemize}
\item whenever $A(x)\in\bigcup\limits_{1\leq i\leq I_{b,\geq}}A(Z^{b,\geq}_i)$, the $p_{A(x)}$-probability for the test to output \up\ is $\geq1-\delta$, and
\item whenever $A(x)\in\bigcup\limits_{1\leq i\leq I_{a,\leq}}A(Z^{a,\leq}_i)$, the $p_{A(x)}$-probability for the test to output \dw\ is $\geq1-\delta$.
\end{itemize}
Situation with a $\delta$-good (right) segment $[a,b]$ is completely similar.
\subsubsection{Bisection estimate: construction}\label{sec:bconstr}
The control parameters of the Bisection  estimate are
\begin{enumerate}
\item positive integer $L$ -- the maximum allowed number of bisection steps,
\item tolerances $\delta\in(0,1)$ and $\varkappa>0$.
\end{enumerate}
The estimate of $f(x)$ ($x$ is the signal underlying our observations: $\omega_t\sim p_{A(x)}$) is given by the following recurrence run on the observation ${\omega}^K=({\omega}_1,...,{\omega}_K)$ which we have at our disposal:
\begin{enumerate}
\item {\bf Initialization.} We {suppose that} a valid upper bound $b_0$ on $\max_{u\in X} f(u)$ and a valid lower bound $a_0$ on $\min_{u\in X}f(u)$ ] {are available;} we assume w.l.o.g. that $a_0<b_0$, otherwise the estimation is trivial. We set $\Delta_0=[a_0,b_0]$ (note that $f(a)\in\Delta_0$).
\item {\bf Bisection Step $\ell$, $1\leq\ell\leq L$}. Given {\sl localizer} $\Delta_{\ell-1}=[a_{\ell-1},b_{\ell-1}]$ with $a_{\ell-1}<b_{\ell-1}$, we act as follows:
\begin{enumerate}
\item\label{step0} Set $c_\ell={1\over 2}[a_{\ell-1}+b_{\ell-1}]$. \par If $c_\ell$ is not upper-feasible, we set $\Delta_\ell=[a_{\ell-1},c_\ell]$ and pass to \ref{lastline}, and if $c_\ell$ is not lower-feasible, we set $\Delta_\ell=[c_\ell,b_{\ell-1}]$ and pass to \ref{lastline}.
    \\
    \underline{Note:} When the rule requires to pass to \ref{lastline}, the set $\Delta_{\ell}\backslash\Delta_{\ell-1}$ does not intersect with $f(X)$; in particular, in this case $f(x)\in\Delta_{\ell}$ provided that $f(x)\in\Delta_{\ell-1}$.
\item\label{step1r} When $c_\ell$ is both upper- and lower-feasible, we check whether the segment $[c_{\ell},b_{\ell-1}]$ is $\delta$-good (right). If it is not the case, we terminate and claim that $f(x)\in\bar{\Delta}:=\Delta_{\ell-1}$, otherwise find $v_\ell$, $c_\ell< v_\ell\leq b_{\ell-1}$, such that the segment $\Delta_{\ell\rg}=[c_{\ell},v_\ell]$ is $\delta$-good (right) $\varkappa$-maximal.
    \\
    \underline{Note:} In terms of the outline of our strategy presented in Section \ref{outline}, termination when the segment $[c_{\ell},b_{\ell-1}]$ is not $\delta$-good (right) corresponds to the case where the current localizer is too small to allow for a separator wide enough to ensure low-risk decision on the left and the right hypotheses.\par
    To find $v_\ell$, we check the candidates  with $v_\ell^k=b_{\ell-1}-k\varkappa$, $k=0,1,...$ until arriving for the first time at segment $[c_{\ell},v_\ell^k]$ which is not $\delta$-good (right), and take, as $v_\ell$, the quantity $v^{k-1}$
    (the resulting value of $v_\ell$ is well defined and clearly meets the above requirements as we clearly have $k\geq1$).
 \item\label{step1l} Similarly, we check whether the segment $[a_{\ell-1},c_\ell]$ is $\delta$-good (left). If it is not the case, we terminate and claim that $f(x)\in\bar{\Delta}:=\Delta_{\ell-1}$, otherwise we find $u_\ell$, $a_{\ell-1}\leq u_\ell< c_\ell$, such that the segment $\Delta_{\ell\lf}=[u_{\ell},c_\ell]$ is $\delta$-good (left) $\varkappa$-maximal.\\
     \underline{Note:} The rules for building $u_\ell$ are completely similar to those for $v_\ell$.
 \item\label{mainstep} We compute $\cT^K_{\Delta_{\ell\rg}\r}({\omega}^K)$ and $\cT^K_{\Delta_{\ell\lf}\l}({\omega}^K)$.
 If $\cT^K_{\Delta_{\ell\rg}\r}({\omega}^K)=\cT^K_{\Delta_{\ell\lf}\l}({\omega}^K)$ (``consensus''), we set
 \begin{equation}\label{NDeltaell}
\Delta_\ell=[a_{\ell},b_{\ell}]=\left\{\begin{array}{ll}[c_\ell,b_{\ell-1}],&\cT^K_{\Delta_{\ell\rg}\r}({\omega}^K)=\up,\\
{[a_{\ell-1},c_\ell]},&\cT^K_{\Delta_{\ell\rg}\r}({\omega}^K)=\dw\\
\end{array}\right.
\end{equation}
and pass to \ref{lastline}. Otherwise (``disagreement'') we terminate and claim that $f(x)\in\bar{\Delta}=[u_\ell,v_\ell]$.
 \item\label{lastline} 
 When $\ell<L$, we pass to step $\ell+1$, otherwise we terminate and claim that $f(x)\in\bar{\Delta}:=\Delta_L$.
 \end{enumerate}
 \item {\bf Output} of the estimation procedure is the segment $\bar{\Delta}$ built upon termination and claimed to contain $f(x)$, see rules
  \ref{step1r} -- \ref{lastline}; the midpoint of this segment is the estimate of $f(x)$ yielded by our procedure.
 \end{enumerate}
\subsubsection{Bisection estimate: Main result}
 \index{Bisection estimate!near-optimality of}
  \begin{proposition}\label{propBisection} Consider the situation described in the beginning of Section \ref{sect:setting}, and let $\epsilon\in(0,1/2)$ be given. Then
 \par\
 {\rm (i) [reliability]} for every positive integer $L$ and every $\varkappa>0$, Bisection with control parameters $L$,
 $
 \delta={\epsilon\over 2L},
 $
 and $\varkappa$ is $(1-\epsilon)$-reliable: for every $x\in X$, the $p_{A(x)}$-probability of the event
 $$
 f(x)\in\bar{\Delta}
 $$
 ($\bar{\Delta}$ is the output of Bisection as defined above) is at least $1-\epsilon$.\par
 {\rm (ii) [near-optimality]} Let $\bar\rho>0$ and positive integer $\bar{K}$ be such that  there exists a $(\bar\rho,\epsilon)$-reliable estimate  $\widehat{f}(\cdot)$ of $f(x)$, $x\in X:=\bigcup_{i\leq I}X_i$, via stationary $\bar{K}$-repeated observation $\omega^{\bar{K}}$ with $\omega_k\sim p_{A(x)}$, $1\leq k\leq \bar{K}$. Given $\rho>2\bar\rho$, the Bisection  estimate utilizing stationary $K$-repeated observations, with
\begin{equation}\label{Kselectedas}
K=\left\rfloor{2\ln(2LNI/\epsilon)\over \ln([4\epsilon(1-\epsilon)]^{-1})}\bar{K}\right\lfloor,
\end{equation}
the control parameters of the estimate being
\begin{equation}\label{Nsetup}
L=\left\rfloor {\log}_2\left({b_0-a_0\over2\rho}\right)\right\lfloor,\;\;\delta={\epsilon\over 2L},\;\;\varkappa=\rho-2\bar\rho,
\end{equation}
is $(\rho,\epsilon)$-reliable.
\end{proposition}
For proof, see Section \ref{propBisectionproof}.
\par
Note
 that the running time $K$ of Bisection estimate as given by (\ref{Kselectedas}) is just by (at most) logarithmic in $N$, $I$, $L$ and $\epsilon^{-1}$ factor larger than $\bar{K}$, and that $L$ is just logarithmic in $1/\bar\rho$. Assume, for instance, that for some $\gamma>0$  there exist $(\epsilon^\gamma,\epsilon)$ reliable estimates, parameterized by  $\epsilon\in(0,1/2)$, with $\bar{K}=\bar{K}(\epsilon)$. Then Bisection with the volume of observation and control parameters given by (\ref{Kselectedas}),  (\ref{Nsetup}), where  $\rho=3\bar\rho=3\epsilon^\gamma$, and $\bar{K}=\bar{K}(\epsilon)$, is  $(3\epsilon^\gamma,\epsilon)$-reliable and requires  $K=K(\epsilon)$-repeated observations with $\overline{\lim}_{\epsilon\to+0}K(\epsilon)/\bar{K}(\epsilon)\leq2$.
 \subsection{Illustration: {estimating} survival rate}
 {Let $\xi\in \bR_+$ be a random variable representing lifetime. Suppose that our objective is, given $K$ independent {\em indirect} observations of $\xi$ and a value $\tau\in \bR$, estimate the corresponding hazard rate $s_\tau=f_\xi(\tau)/(1-F_\xi(\tau))$ where $f_\xi$ and $F_\xi$ are, respectively, density and cumulative distribution function of $\xi$. Suppose that the density $f_\xi$ is smooth with bounded second derivative, and that observations are subjected to ``mixed'' {\em multiplicative censoring} (see, e.g. \cite{vardi1989,andersen2001,brunel2016,belomestny2017}): the exact value of $\xi_k$ is observed with probability $0\leq \theta\leq 1$, and with complementary probability, the available observation is $\eta_k\xi_k$, where $\eta_k$ is uniformly distributed over $[0,1]$. \par
We assume that after an appropriate discretization, the estimation problem can be reformulated as follows: let $x$ be the distribution of the (discrete-valued) lifetime taking values in $S=\{1,2,...,M\}$. We define the corresponding hazard rate $s_{j}[x]$ (the conditional probability of the lifetime to be exactly $j$ given that it is at least $j$) according to
\[
s_{j}(x)={x_j\over\sum_{i=j}^M x_i},\,1\leq j\leq M.
\]
Our objective is to estimate $s_j[x]$, given $K$ independent observations $\omega_k$ with distribution $\mu=Ax$, where $A\in \bR^{M\times M}$ is a given column-stochastic matrix.
\par
We use the following setup:
 \begin{itemize}
 \item $X=\{x\in\bR^m:\;x_i\geq(3M)^{-1};\,\sum_{i=1}^Mx_i=1;\,|x_{i-1}-2x_i+x_{i+1}|\leq 2M^{-2},1<i<M\}$;
 \item $A=\theta I_M+(1-\theta)R$, where $R$ is upper-triangular matrix with the $i$-th column $(\underbrace{i^{-1},...,i^{-1}}_{i},0,...,0)^T$.
 \end{itemize}
 }
 {
For various combinations of $\theta$ and $K$ we carried out 100 simulations of bisection estimation. In each simulation, we first selected $x\in X$ at random, drew $K$ observations $\omega_t$, $t=1,...,K$, from the distribution $Ax$, and then ran Bisection on these observations.
Plots in Figure \ref{fig:hr_tK} illustrate some typical results of our experiments.}
\begin{figure}[h!]
\begin{tabular}{cc}
\includegraphics[scale=0.55]{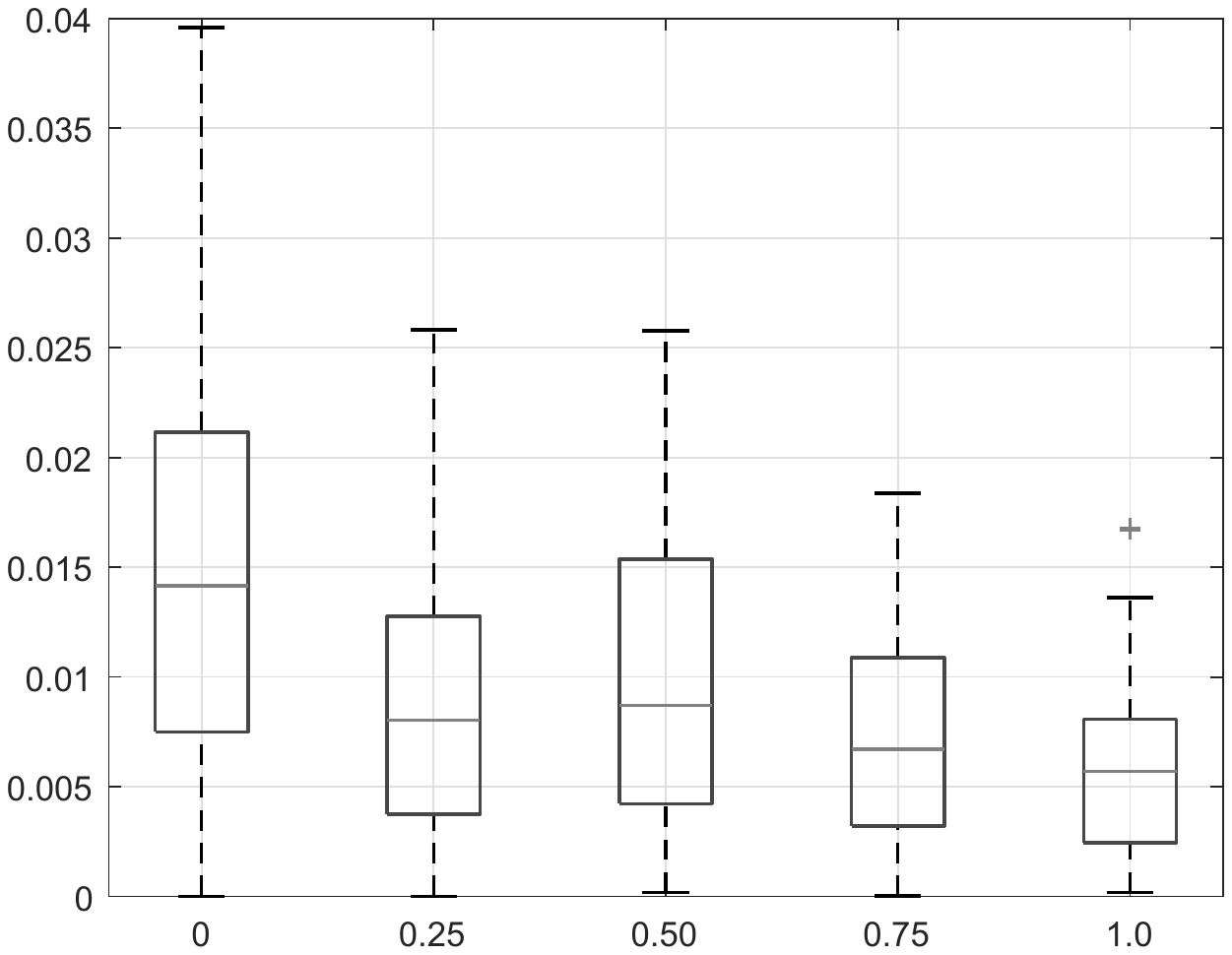}&\includegraphics[scale=0.55]{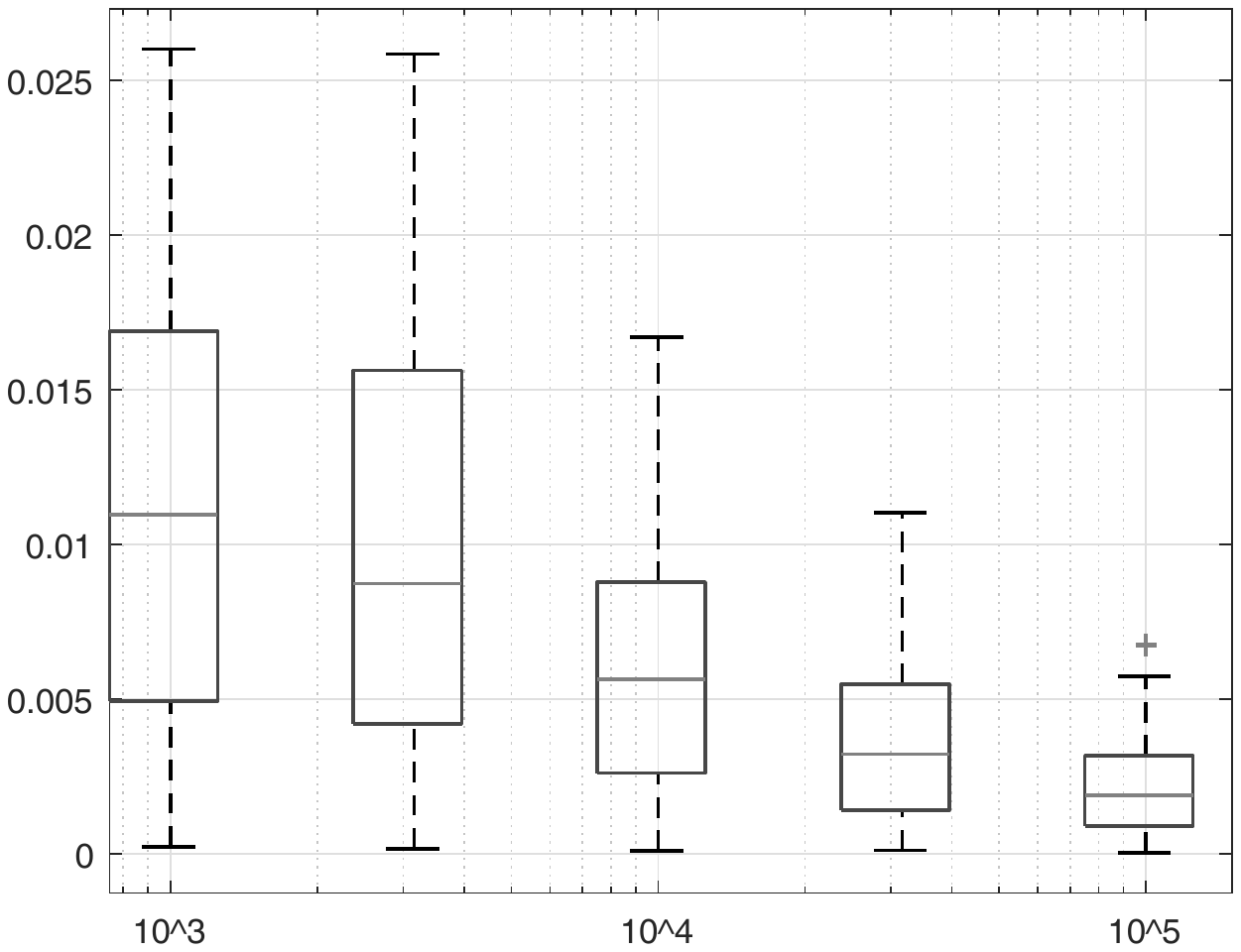}\\
{\small $(a)$}&{\small $(b)$}
\end{tabular}
\caption{\small Boxplot of emprirical error distribution of Bisection estimate over 100 random estimation problems. (a) For $K=10\,000$, hasard rate estimation error as function of $\theta\in \{0,0.25,0.5,0.75,1\}$; (b) estimation error as function of $K$ for $\theta=0.9$. In these experiments, the initial risk -- the half-width of the initial localizer -- is equal to $0.0524$.}
\label{fig:hr_tK}
\end{figure}
}


\appendix
\section{Proofs}
\subsection{Proof of Proposition \ref{proplin}}\label{pr:00}
{\bf Proof.} Let
the common distribution $p$ of independent across $k$ components $\omega_k$ of $\omega^K$ be  $p_{A_\ell(u)}$ for some $\ell\leq I$ and $u\in X_\ell$. Let us fix these $\ell$ and $u$, let $\mu=A_\ell(u)$, and let $p^K$ stand for the distribution of $\omega^K$.
\paragraph{1$^0$.} We have
\[
\begin{array}{rcl}
\Psi_{\ell,+}(\alpha_{\ell j},\phi_{\ell j})&=&\max_{x\in X_\ell}\left[K\alpha_{\ell j} \Phi_\cO(\phi_{\ell j}/\alpha_{\ell j},A_\ell(x))-g^Tx\right]+\alpha_{\ell j}\ln(2I/\epsilon)\\
&\geq&  K\alpha_{\ell j} \Phi_\cO(\phi_{\ell j}/\alpha_{\ell j},\mu)-g^Tu+\alpha_{\ell j}\ln(2I/\epsilon){\hbox{\ [since $u\in X_\ell$ and $\mu=A_\ell(u)$]}}\\
&=&K\alpha_{\ell j}\ln\left(\int \exp\{\phi_{\ell j}(\omega)/\alpha_{\ell j}\}p_\mu(\omega)P(d\omega)\right)-g^Tu+\alpha_{\ell j}\ln(2I/\epsilon)\hbox{\ [definition of $\Phi_\cO$]}\\
&=&\alpha_{\ell j} \ln\left(\bE_{\omega^K\sim p^K}\left\{\exp\{\alpha_{\ell j}^{-1}\sum_k\phi_{\ell j}(\omega_k)\}\right\}\right)-g^Tu+\alpha_{\ell j}\ln(2I/\epsilon)\\
&=&\alpha_{\ell j} \ln\left(\bE_{\omega^K\sim p^K}\left\{\exp\{\alpha_{\ell j}^{-1}[g_{\ell j}(\omega^K)-\varkappa_{\ell j}]\}\right\}\right)-g^Tu+\alpha_{\ell j}\ln(2I/\epsilon)\\
&=&\alpha_{\ell j} \ln\left(\bE_{\omega^K\sim p^K}\left\{\exp\{\alpha_{\ell j}^{-1}[g_{\ell j}(\omega^K)-g^Tu-\rho_{\ell j}]\}\right\}\right) +\rho_{\ell j}-\varkappa_{\ell j}+\alpha_{\ell j}\ln(2I/\epsilon)\\
&\geq& \alpha_{\ell j} \ln\left(\Prob_{\omega^K\sim p^K}\left\{g_{\ell j}(\omega^K)> g^Tu+\rho_{\ell j}\right\}\right) +\rho_{\ell j}-\varkappa_{\ell j}+\alpha_{\ell j}\ln(2I/\epsilon)
\end{array}
\]
so that
\[
\begin{array}{rcl}
\alpha_{\ell j}\ln\left(\Prob_{\omega^K\sim p^K}\left\{g_{\ell j}(\omega^K)> g^Tu+\rho_{\ell j}\right\}\right)&\leq&
\Psi_{\ell,+}(\alpha_{\ell j},\phi_{\ell j})+\varkappa_{\ell j}-\rho_{\ell j}+\alpha_{\ell j}\ln({\epsilon\over 2I})\\&=&\alpha_{\ell j}\ln({\epsilon\over 2I}) \hbox{\ [by (\ref{gij})]},
\end{array}
\]
and we arrive at
\begin{equation}\label{upperbound}
\Prob_{\omega^K\sim p^K}\left\{g_{\ell j}(\omega^K)> \rho_{\ell j}+g^Tu\right\}\leq {\epsilon\over 2I}.
\end{equation}
Similarly,
\[
\begin{array}{rcl}
\Psi_{\ell,-}(\alpha_{i \ell},\phi_{i\ell})&=&\max_{y\in X_\ell}\left[K\alpha_{i\ell} \Phi_\cO(-\phi_{i\ell}/\alpha_{i\ell},A_\ell(y))+g^Ty\right]+\alpha_{i\ell }\ln(2I/\epsilon)\\
&\geq&  K\alpha_{i\ell} \Phi_\cO(-\phi_{i\ell}/\alpha_{i\ell},\mu)+g^Tu+\alpha_{i\ell}\ln(2I/\epsilon){\hbox{\ [since $u\in X_\ell$ and $\mu=A_\ell(u)$]}}\\
&=&K\alpha_{i\ell}\ln\left(\int \exp\{-\phi_{i\ell}(\omega)/\alpha_{i\ell}\}p_\mu(\omega)P(d\omega)\right)+g^Tu+\alpha_{i\ell}\ln(2I/\epsilon)\hbox{\ [definition of $\Phi_\cO$]}\\
&=&\alpha_{i\ell} \ln\left(\bE_{\omega^K\sim p^K}\left\{\exp\{-\alpha_{i\ell}^{-1}\sum_k\phi_{i\ell}(\omega_k)\}\right\}\right)+g^Tu+\alpha_{i\ell}\ln(2I/\epsilon)\\
&=&\alpha_{i\ell} \ln\left(\bE_{\omega^K\sim p^K}\left\{\exp\{\alpha_{i\ell}^{-1}[-g_{i\ell}(\omega^K)+\varkappa_{i\ell}]\}\right\}\right)+g^Tu+\alpha_{i\ell}\ln(2I/\epsilon)\\
&=&\alpha_{i\ell} \ln\left(\bE_{\omega^K\sim p^K}\left\{\exp\{\alpha_{i\ell}^{-1}[-g_{i\ell}(\omega^K)+g^Tu-\rho_{i\ell}]\}\right\}\right) +\rho_{i\ell}+\varkappa_{i\ell}+\alpha_{i\ell}\ln(2I/\epsilon)\\
&\geq& \alpha_{i\ell} \ln\left(\Prob_{\omega^K\sim p^K}\left\{g_{i\ell}(\omega^K)< g^Tu-\rho_{i\ell}\right\}\right) +\rho_{i\ell}+\varkappa_{i\ell}+\alpha_{i\ell}\ln(2I/\epsilon),\end{array}\]
implying that
\[\begin{array}{rcl}
\alpha_{i\ell}\ln\left(\Prob_{\omega^K\sim p^K}\left\{g_{i\ell}(\omega^K)< g^Tu-\rho_{i\ell}\right\}\right)&\leq&
\Psi_{\ell,-}(\alpha_{i\ell},\phi_{i\ell})-\varkappa_{i\ell}-\rho_{i\ell}+\alpha_{i\ell}\ln({\epsilon\over 2I})\\
&=&\alpha_{i\ell}\ln({\epsilon\over 2I}) \hbox{\ [by (\ref{gij})]},
\end{array}
\]
and we conclude that
\begin{equation}\label{lowerbound}
\Prob_{\omega^K\sim p^K}\left\{g_{i\ell}(\omega^K)<g^Tu-\rho_{i\ell}\right\}\leq {\epsilon\over 2I}.
\end{equation}
\paragraph{2$^0$.} Let
$$
\cE=\{\omega^K: g_{\ell j}(\omega^K) \leq g^Tu+\rho_{\ell j}, \;g_{i\ell}(\omega^K) \geq g^Tu-\rho_{i\ell},\; 1\leq i,j\leq I\}.
$$
From (\ref{upperbound}), (\ref{lowerbound}) and the union bound it follows that $p^K$-probability of the event $\cE$ is $\geq1-\epsilon$. As a result, all we need to complete the proof of Proposition is to verify that for all $\omega^K\in \cE$,
\begin{equation}\label{verifythat}
|\widehat{g}(\omega^K)-g^Tu|\leq \rho_\ell.
\end{equation}
Indeed, let us fix $\omega^K \in\cE$, and let $E$ be the $I\times I$ matrix with entries $E_{ij}=g_{ij}(\omega^K)$, $1\leq i,j\leq I$. The quantity $r_i$, see (\ref{rcij}), is the maximum of entries in $i$-th row of $E$, and the quantity $c_j$ is the minimum of entries in $j$-th column of $E$. In particular, $r_i\geq E_{ij}\geq c_j$ for all $i,j$, implying that {$r_i\geq c_\ell$ and $c_j\leq r_\ell$} for all $i,j$.
Now, since $\omega^K\in \cE$, we have { for all $j$:}
\[
E_{\ell j}=g_{\ell j}(\omega^K)\leq g^Tu+\rho_{\ell j}\leq g^Tu+\rho_\ell,
\] implying that $r_\ell=\max_j E_{\ell j}{\leq g^Tu+\rho_\ell}$.
Similarly, $\omega^K\in\cE$ implies that  for all $i$
\[
E_{i\ell}=g_{i\ell}(\omega^K) \geq g^Tu-\rho_{i\ell}\geq g^Tu-\rho_\ell,
\] so that $c_\ell=\min_iE_{i\ell}{\geq g^Tu-\rho_\ell}$.  {
 We have $r_*:=\min_{i}r_i\leq r_\ell$, and, as we have already seen, $r_*\geq c_\ell$, implying  that $r_*$ belongs to $\Delta_\ell=[g^Tu-\rho_\ell,g^Tu+\rho_\ell]$. By similar argument, $c_*:=\max_jc_j\in\Delta_\ell$ as well.}
Finally, $\widehat{g}(\omega^K)={1\over 2}[r_*+c_*]$, that is, $\widehat{g}(\omega^K)\in\Delta_\ell$, and \rf{verifythat} follows. \qed

\subsection{Proof of Proposition \ref{propnearoptlin}}\label{pr:01}
{\bf 1$^0$.} Observe that $\Opt_{ij}(K)$ is the saddle point value in the convex-concave saddle point problem:
{\small $$
\Opt_{ij}(K)=\inf_{\alpha>0,\phi\in\cF}\max_{x\in X_i,y\in X_j}\left[{\half}K\alpha\left\{\Phi_\cO(\phi/\alpha;A_i(x))+\Phi_\cO(-\phi/\alpha;A_j(y))\right\}+\half g^T[y-x]+\alpha\ln(2I/\epsilon)\right].
$$}\noindent
The domain of the maximization variable is compact and the cost function is continuous on its domain, whence, by Sion-Kakutani Theorem, we have also
\be\begin{array}{rcl}
\Opt_{ij}(K)&=&\max\limits_{x\in X_i,y\in X_j}\Theta_{ij}(x,y),\\
\Theta_{ij}(x,y)&=&\inf\limits_{\alpha>0,\phi\in\cF}\left[\half K\alpha\left\{
\Phi_\cO(\phi/\alpha;A_i(x))+\Phi_\cO(-\phi/\alpha;A_j(y))\right\}+\alpha\ln(2I/\epsilon)\right]
+\half g^T[y-x].
\end{array}\ee{wehavealso}
We have
$$
\begin{array}{rcl}
\Theta_{ij}(x,y)
&=&\inf\limits_{\alpha>0,\psi\in\cF}\left[
{\half}K\alpha\left\{
\Phi_\cO(\psi;A_i(x))+\Phi_\cO(-\psi;A_j(y))\right\}+\alpha\ln(2I/\epsilon)\right]+\half g^T[y-x]\\
&=&\inf\limits_{\alpha>0}\left[
{\half}\alpha K\inf\limits_{\psi\in \cF}\left\{
\Phi_\cO(\psi;A_i(x))+\Phi_\cO(-\psi;A_j(y))\right\}+\alpha\ln(2I/\epsilon)\right]+{\half}g^T[y-x]
\end{array}
$$
Given $x\in X_i$, $y\in X_j$ and setting $\mu=A_i(x)$, $\nu=A_j(y)$, we obtain
\bse
\inf_{\psi\in \cF}[
\Phi_\cO(\psi;A_i(x))+\Phi_\cO(-\psi;A_j(y))]
&=&\inf_{\psi\in\cF}\left[\ln\left(\int\exp\{\psi(\omega)\}p_\mu(\omega)P(d\omega)\right)\right.\\&&+\left.
\ln\left(\int\exp\{-\psi(\omega)\}p_\nu(\omega)P(d\omega)\right)\right].
\ese Since $\cO$ is a good o.s., the function $\bar{\psi}(\omega)={1\over 2}\ln(p_\nu(\omega)/p_\mu(\omega))$ belongs to $\cF$, and
\bse
\lefteqn{\inf_{\psi\in\cF}\left[\ln\left(\int\exp\{\psi(\omega)\}p_\mu(\omega)P(d\omega)\right)+
\ln\left(\int\exp\{-\psi(\omega)\}p_\nu(\omega)P(d\omega)\right)\right]}\\
&=&\inf_{\delta\in\cF}\left[\ln\left(\int\exp\{\bar{\psi}(\omega)+\delta(\omega)\}p_\mu(\omega)P(d\omega)\right)+
\ln\left(\int\exp\{-\bar{\psi}(\omega)-\delta(\omega)\}p_\nu(\omega)P(d\omega)\right)\right]\\
&=&\inf_{\delta\in\cF} \underbrace{\left[\ln\left(\int\exp\{\delta(\omega)\}\sqrt{p_\mu(\omega)p_\nu(\omega)}P(d\omega)\right)+
\ln\left(\int\exp\{-\delta(\omega)\}\sqrt{p_\mu(\omega)p_\nu(\omega)}P(d\omega)\right)\right]}_{f(\delta)}.
\ese
Observe that $f(\delta)$ clearly is a convex and even function of $\delta\in\cF$; as such, it attains its minimum over $\delta\in\cF$ when $\delta=0$. The bottom line is that
\begin{equation}\label{bottomlineis}
\inf_{\psi\in \cF}[
\Phi_\cO(\psi;A_i(x))+\Phi_\cO(-\psi;A_j(y))]=2\ln\left(\int\sqrt{p_{A_i(x)}(\omega)p_{A_j(y)}(\omega)}P(d\omega)\right),
\end{equation}
and
\bse\Theta_{ij}(x,y)&=&\inf_{\alpha>0}\alpha\left[K\ln\left(\int\sqrt{p_{A_i(x)}(\omega)p_{A_j(y)}(\omega)}P(d\omega)\right)+\ln(2I/\epsilon)\right]+{\half}g^T[y-x]\\
&=&\left\{\begin{array}{ll}{\half}g^T[y-x]&,K\ln\left(\int\sqrt{p_{A_i(x)}(\omega)p_{A_j(y)}(\omega)}P(d\omega)\right)+\ln(2I/\epsilon)\geq0,\\
-\infty&,\hbox{otherwise}.\\
\end{array}\right.
\ese
This combines with (\ref{wehavealso}) to imply that
\begin{equation}\label{wehaveeven}
\Opt_{ij}(K)=\max_{x,y}\left\{\half g^T[y-x]: x\in X_i,y\in X_j, \left[\int\sqrt{p_{A_i(x)}(\omega)p_{A_j(y)}(\omega)}P(d\omega)\right]^K\geq {\epsilon\over 2I}\right\}.
\end{equation}
\paragraph{2$^0$.} We claim that under the premise of Proposition, for all $i,j$, $1\leq i,j\leq I$, one has
\[
\Opt_{ij}(K)\leq \Risk^*_\epsilon(\bar{K}),
\]
implying the validity of (\ref{onehasone}).
Indeed, assume that for some pair $i,j$ the opposite inequality holds true:
$$
\Opt_{ij}(K)>\Risk^*_\epsilon(\bar{K}),
$$
and let us lead this assumption to a contradiction. Under our assumption optimization problem in (\ref{wehaveeven}) has a feasible solution $(\bar{x},\bar{y})$ such that
\begin{equation}\label{suchthat}
r:={\half}g^T[\bar{y}-\bar{x}] >\Risk^*_\epsilon(\bar{K}),
\end{equation}
implying, due to the origin of $\Risk^*_\epsilon(\bar{K})$, that there exists an estimate $\widehat{g}(\omega^{\bar{K}})$ such that for $\mu=A_i(\bar{x})$, $\nu=A_j(\bar{y})$ it holds
$$
\begin{array}{rcl}
\Prob_{\omega^{\bar{K}}\sim p_\nu^{\bar{K}}}\left\{\widehat{g}(\omega^{\bar{K}})\leq{\half}g^T[\bar{x}+\bar{y}]\right\}&\leq&\Prob_{\omega^{\bar{K}}\sim p_\nu^{\bar{K}}}\left\{|\widehat{g}(\omega^{\bar{K}})-g^T\bar{y}|\geq r\right\}\leq\epsilon\\
\Prob_{\omega^{\bar{K}}\sim p_\mu^{\bar{K}}}\left\{\widehat{g}(\omega^{\bar{K}})\geq{\half}g^T[\bar{x}+\bar{y}]\right\}&\leq&\Prob_{\omega^{\bar{K}}\sim p_\mu^{\bar{K}}}\left\{|\widehat{g}(\omega^{\bar{K}})-g^T\bar{x}|\geq r\right\}\leq\epsilon,
\end{array}
$$
so that we can decide on two simple hypotheses stating that observation $\omega^{\bar{K}}$ obeys distribution $p_\mu^{\bar{K}}$, resp.,  $p_\nu^{\bar{K}}$, with risk $\leq\epsilon$. Therefore,
$$
\int\min\left[p_\mu^{\bar{K}}(\omega^{\bar{K}}),p_\nu^{\bar{K}}(\omega^{\bar{K}})\right]P^{\bar{K}}(d\omega^{\bar{K}})\leq2\epsilon.
\eqno{[P^{\bar{K}}=\underbrace{P\times...\times P}_{\bar{K}}]}
$$
Hence, when setting $p^{\bar{K}}_\theta(\omega^{\bar{K}})=\prod_kp_\theta(\omega_k)$, we have
\[
\begin{array}{rcl}
\multicolumn{3}{l}{\left[\int\sqrt{p_\mu(\omega)p_\nu(\omega)}P(d\omega)\right]^{\bar{K}}=
\int\sqrt{p_\mu^{\bar{K}}(\omega^{\bar{K}})p_\nu^{\bar{K}}(\omega^{\bar{K}})}P^{\bar{K}}(d\omega^{\bar{K}})}\\
&=&\int\sqrt{\min\left[p_\mu^{\bar{K}}(\omega^{\bar{K}}),p_\nu^{\bar{K}}(\omega^{\bar{K}})\right]}\sqrt{\max\left[p_\mu^{\bar{K}}(\omega^{\bar{K}}),
p_\nu^{\bar{K}}(\omega^{\bar{K}})\right]}P^{\bar{K}}(d\omega^{\bar{K}})\\
&\leq& \left[\int\min\left[p_\mu^{\bar{K}}(\omega^{\bar{K}}),p_\nu^{\bar{K}}(\omega^{\bar{K}})\right]P^{\bar{K}}(d\omega^{\bar{K}})\right]^{1/2}
\left[\int\max\left[p_\mu^{\bar{K}}(\omega^{\bar{K}}),p_\nu^{\bar{K}}(\omega^{\bar{K}})\right]P^{\bar{K}}(d\omega^{\bar{K}})\right]^{1/2}\\
&=&\left[\int\min\left[p_\mu^{\bar{K}}(\omega^{\bar{K}}),p_\nu^{\bar{K}}(\omega^{\bar{K}})\right]P^{\bar{K}}(d\omega^{\bar{K}})\right]^{1/2}
\\&&\times\left[\int\left[p_\mu^{\bar{K}}(\omega^{\bar{K}})+p_\nu^{\bar{K}}(\omega^{\bar{K}})-
\min\left[p_\mu^{\bar{K}}(\omega^{\bar{K}}),p_\nu^{\bar{K}}(\omega^{\bar{K}})\right]\right]P^{\bar{K}}(d\omega^{\bar{K}})\right]^{1/2}\\
&=&\left[\int\min\left[p_\mu^{\bar{K}}(\omega^{\bar{K}}),p_\nu^{\bar{K}}(\omega^{\bar{K}})\right]P^{\bar{K}}(d\omega^{\bar{K}})\right]^{1/2}\left[2-
\int\min\left[p_\mu^{\bar{K}}(\omega^{\bar{K}}),p_\nu^{\bar{K}}(\omega^{\bar{K}})\right]P^{\bar{K}}(d\omega^{\bar{K}})\right]^{1/2}\\
&\leq& 2\sqrt{\epsilon(1-\epsilon)}.
\end{array}
\]
Consequently,
$$
\left[\int\sqrt{p_\mu(\omega)p_\nu(\omega)}P(d\omega)\right]^{K}\leq [2\sqrt{\epsilon(1-\epsilon)}]^{K/\bar{K}}<{\epsilon\over 2I},
$$
which is the desired contradiction (recall that $\mu=A_i(\bar{x})$, $\nu=A_j(\bar{y})$ and $(\bar{x},\bar{y})$ is feasible for  (\ref{wehaveeven})).
\paragraph{3$^0$.} Now let us prove that under the premise of Proposition, (\ref{onehastwo}) takes place. To this end let us set
\begin{equation}\label{Gammaijis}
w_{ij}(s)=\max_{x\in X_j,y\in X_j}\bigg\{{\half}g^T[y-x]: \bar{K}\underbrace{\ln\left(\int\sqrt{p_{A_i(x)}(\omega)p_{A_j(y)}(\omega)}P(d\omega)\right)}_{H(x,y)}+s\geq0\bigg\}.
\end{equation}
As we have seen in item 1$^0$, see (\ref{bottomlineis}), one has
$$
H(x,y)=\inf_{\psi\in \cF}{\half}\left[\Phi_\cO(\psi;A_i(x))+\Phi_\cO(-\psi,A_j(y))\right],
$$
that is, $H(x,y)$ is the infimum of a parametric family of concave functions of $(x,y)\in X_i\times X_j$ and as such is concave. Besides this,
the optimization problem in (\ref{Gammaijis}) is feasible whenever $s\geq0$, a feasible solution being $y=x=x_{ij}$. At this feasible solution we have $g^T[y-x]=0$, implying that  $w_{ij}(s)\geq0$ for $s\geq0$. Observe also that from concavity of  $H(x,y)$ it follows that $w_{ij}(s)$ is concave on the ray $\{s\geq0\}$. Finally, we claim that
\begin{equation}\label{weclaim}
w_{ij}(\bar{s})\leq \Risk^*_\epsilon(\bar{K}),\,\bar{s}=-\ln(2\sqrt{\epsilon(1-\epsilon)}).
\end{equation}
Indeed, $w_{ij}(s)$ is nonnegative, concave and bounded (since $X_i,X_j$ are compact) on $\bR_+$, implying that $w_{ij}(s)$ is continuous on $\{s>0\}$. Assuming, on the contrary to our claim, that $w_{ij}(\bar{s})>\Risk^*_\epsilon(\bar{K})$, there exists $s'\in(0,\bar{s})$ such that $w_{ij}(s')>\Risk^*_\epsilon(\bar{K})$ and thus there exist
$\bar{x}\in X_i$, $\bar{y}\in X_j$ such that $(\bar{x},\bar{y})$ is feasible for the optimization problem specifying $w_{ij}(s')$ and (\ref{suchthat}) takes place. We have seen in item 2$^0$ that the latter relation implies that for $\mu=A_i(\bar{x})$, $\nu=A_j(\bar{y})$ it holds
$$
\left[\int\sqrt{p_\mu(\omega)p_\nu(\omega)}P(d\omega)\right]^{\bar{K}}\leq 2\sqrt{\epsilon(1-\epsilon)},
$$
that is,
$$
\bar{K}\ln\left(\int\sqrt{p_\mu(\omega)p_\nu(\omega)}P(d\omega)\right)+\bar{s}\leq 0,
$$
whence
$$
\bar{K}\ln\left(\int\sqrt{p_\mu(\omega)p_\nu(\omega)}P(d\omega)\right)+s'< 0,
$$
contradicting the fact that $(\bar{x},\bar{y})$ is feasible for the optimization problem specifying $w_{ij}(s')$.
\par
It remains to note that (\ref{weclaim}) combines with concavity of $w_{ij}(\cdot)$ and the relation $w_{ij}(0)\geq0$ to imply that
$$
w_{ij}(\ln(2I/\epsilon))\leq \vartheta w_{ij}(\bar{s})\leq\vartheta\Risk^*_\epsilon(\bar{K}),\;\;\vartheta=\ln(2I/\epsilon)/\bar{s}
={2\ln(2I/\epsilon)\over \ln([4\epsilon(1-\epsilon)]^{-1})}.
$$
Invoking (\ref{wehaveeven}), we conclude that
$$
\Opt_{ij}(\bar{K})=w_{ij}(\ln(2I/\epsilon))\leq \vartheta\Risk^*_\epsilon(\bar{K})\,\forall i,j.
$$
 Finally, from (\ref{wehaveeven}) it immediately follows that $\Opt_{ij}(K)$ is nonincreasing in $K$ (since as $K$ grows, the feasible set of the right hand side optimization problem in (\ref{wehaveeven}) shrinks), that is,
$$
K\geq \bar{K}\Rightarrow \Opt(K)\leq \Opt(\bar{K})=\max_{i,j}\Opt_{ij}(\bar{K}) \leq \vartheta\Risk^*_\epsilon(\bar{K}),
$$
and (\ref{onehastwo}) follows. \qed
\subsection{Proof of Proposition \ref{propBisection}}\label{propBisectionproof}
\subsubsection{Proof of Proposition \ref{propBisection}(i)}\label{Nanalysis}
 We call step $\ell$ {\sl {constructive}}, if at this step rule \ref{mainstep} is invoked.
\paragraph{1$^0$.}
 Let $x\in X$ be the true signal underlying our observation ${\omega}^K$, so that ${\omega}_1,...,{\omega}_K$ are independently of each other drawn from the distribution $p_{A(x)}$. Consider the ``ideal'' Bisection given by exactly the same rules as the procedure described in Section \ref{sec:bconstr} (in the sequel, we refer to the latter {as to the ``actual''} one), up to the fact that  tests $\cT^K_{\Delta_{\ell\rg}\r}(\cdot)$, $\cT^K_{\Delta_{\ell\lf}\l}(\cdot)$ in rule \ref{mainstep} {are replaced by the rules}
 $$
 T^*_{\Delta_{\ell\rg}\r}=T^*_{\Delta_{\ell\lf}\l}=\left\{\begin{array}{ll}\up,&f(x)>c_\ell\\
 \dw,&f(x)\leq c_\ell\\
 \end{array}
 \right.
 $$
Marking by $^*$ the entities produced by the resulting {\sl deterministic} procedure, we arrive at a sequence of nested segments $\Delta^*_\ell=[a_\ell^*,b_\ell^*]$, $0\leq\ell\leq L^*\leq L$, along with subsegments $\Delta_{\ell\rg}^*=[c_\ell^*,v_\ell^*]$, $\Delta_{\ell\lf}^*=[u_\ell^*,c_\ell^*]$ of $\Delta_{\ell-1}^*$, defined for all $^*$-{constructive} {steps} $\ell$, and the output segment $\bar{\Delta}^*$ claimed to contain $f(x)$. Note that the ideal procedure cannot terminate due to a disagreement, and that $f(x)$, as is immediately seen, is contained in all segments $\Delta_\ell^*$, $0\leq\ell\leq L^*$, same as $f(x)\in\bar{\Delta}^*$.\par
Let $\cL^*$ be the set of all $^*$-{constructive} values of $\ell$. For $\ell\in\cL^*$, let the event
$\cE_\ell[x]$ parameterized by $x$ be defined as follows:
\begin{equation}\label{badevent}
\cE_\ell[x]=\left\{\begin{array}{ll}\{\omega^K: \cT^K_{\Delta^*_{\ell\rg}\r}(\omega^K)=\up\ \hbox{or}\ \cT^K_{\Delta^*_{\ell\lf}\l}(\omega^K)=\up\},&f(x)\leq u_\ell^*\\
\{\omega^K:\cT^K_{\Delta^*_{\ell\rg}\r}(\omega^K)=\up\},&u_\ell^*<f(x)\leq c_\ell^*\\
\{\omega^K:\cT^K_{\Delta^*_{\ell\lf}\l}(\omega^K)=\dw\},&c_\ell^*<f(x)<v_\ell^*\\
\{\omega^K: \cT^K_{\Delta^*_{\ell\rg}\r}(\omega^K)=\dw\ \hbox{or}\ \cT^K_{\Delta^*_{\ell\lf}\l}(\omega^K)=\dw\},&f(x)\geq v_\ell^*\\
\end{array}\right.
\end{equation}
\paragraph{2$^0$.} {Observe that by construction and in view of  Proposition  \ref{basicprop} we have}
\begin{equation}\label{Nwehave}
\forall \ell\in\cL^*: \Prob_{\omega^K\sim p_{A(x)}\times...\times p_{A(x)}}\{\cE_\ell[x]\}\leq 2\delta.
\end{equation}
\begin{quote}
Indeed, let $\ell\in\cL^*$.
\begin{itemize}
\item When $f(x)\leq u_\ell^*$,  we have $x\in X$ and $f(x)\leq u_\ell^*\leq c_\ell^*$, implying that $\cE_\ell[x]$ takes place only when  either the left test $\cT_{\Delta^*_{\ell\lf}\l}^K$, or the right  test
    $\cT_{\Delta^*_{\ell\rg}\r}^K$, or both, did not accept true -- left -- hypotheses from the pairs of right and left hypotheses the tests were applied to. Since the corresponding intervals ($[u_\ell^*,c_\ell^*]$ for the left side test, $[c_\ell^*,v_\ell^*]$ for the right side one) are $\delta$-good left/right, respectively, the risks of the tests do not exceed $\delta$, and the $p_{A(x)}$-probability of the event $\cE_\ell[x]$ is at most $2\delta$;
\item when $u_\ell^*< f(x)\leq c_\ell^*$, the event $\cE_\ell[x]$ takes place only when  the right  test $\cT_{\Delta^*_{\ell\rg}\r}^K$ does not accept true -- left --
{hypothesis}; similarly to the above, this can happen with $p_{A(x)}$-probability at most $\delta$;
\item when $c_\ell <f(x)\leq v_\ell$, the event $\cE_\ell[x]$ takes place only when the left test $\cT_{\Delta^*_{\ell\lf}\l}^K$ does not accept true -- right -- hypothesis, which, again, happens with $p_{A(x)}$-probability $\leq\delta$;
\item finally, when $f(x)>v_\ell$, the event $\cE_{\ell}[x]$ takes place only when either the left test $\cT_{\Delta^*_{\ell\lf}\l}^K$, or the right test
    $\cT_{\Delta^*_{\ell\rg}\r}^K$, or both, does not accept the true -- right -- hypothesis from the pair of right and left hypotheses the test was applied to; same as above, this can happen with $p_{A(x)}$-probability at most $2\delta$.
\end{itemize}
\end{quote}
\paragraph{3$^0$.} Let $\bar{L}=\bar{L}(\bar{\omega}^K)$ be the last step of the ``{actual}'' estimating procedure as run on the observation $\bar{\omega}^K$. We claim that the following holds true:
\begin{lemma}\label{lem:appendix000}
{Let $\cE:=\bigcup_{\ell\in\cL^*}\cE_\ell[x]$, so that the $p_{A(x)}$-probability of the event $\cE$, the observations stemming from $x$, is at most
$$2\delta L=\epsilon$$
 by {\rm (\ref{Nwehave})}. Assume that $\bar{\omega}^K\not\in\cE$. Then $\bar{L}({\omega}^K)\leq L^*$, and just two cases are possible:
\item {\rm (A)}
The actual estimating procedure did not terminate by disagreement. In this case $\bar{L}({\omega}^K)=L^*$, and the trajectories of the ideal and the actual Bisections are identical} (same localizers, same {constructive} steps, same output segments, etc.); in particular, $f(x)\in\Bar{\Delta}$;
\item {\rm (B)} The actual estimating procedure terminated due to a disagreement. Then $\Delta_\ell=\Delta^*_\ell$ for $\ell<\bar{L}$, and $f(x)\in\bar{\Delta}$.
\end{lemma}
In view of {{\rm (A)} and {\rm (B),}} the $p_{A(x)}$-probability of the event $f(x)\in\Bar{\Delta}$ is at least $1-\epsilon$, as claimed in Proposition \ref{propBisection}.

{\bf Proof of the lemma.} Note that the actions at step $\ell$ in the ideal and the actual procedures depend solely on $\Delta_{\ell-1}$ and on the outcome of rule \ref{mainstep}. Taking into account that $\Delta_0=\Delta_0^*$, all we need to verify is the following:
\begin{quote}
(!)
{\sl Let $\bar{\omega}^K\not\in\cE$, and let $\ell\leq L^*$ be such that $\Delta_{\ell-1}=\Delta^*_{\ell-1}$, whence also $u_\ell=u_\ell^*,\;c_\ell=c_\ell^*$ and $v_\ell=v_\ell^*$. Assume that $\ell$ is {constructive} (given that $\Delta_{\ell-1}=\Delta_{\ell-1}^*$, this may happen if and only if $\ell$ is $^*$-{constructive} as well). Then either
\par
-- at step $\ell$ the actual procedure terminates due to disagreement, in which case $f(x)\in\bar{\Delta}$, or\par
-- there was no disagreement at step $\ell$, in which case $\Delta_\ell$ as given by  {\rm (\ref{NDeltaell})} is identical to $\Delta_\ell^*$ as given by the ideal counterpart of {\rm (\ref{NDeltaell})} in the case of $\Delta_{\ell-1}^*=\Delta_{\ell-1}$, that is, by the rule}
\begin{equation}\label{NDeltaellstar}
\Delta_\ell^*=\left\{\begin{array}{ll}[c_\ell,b_{\ell-1}],&f(x)>c_\ell,\\
{[a_{\ell-1},c_\ell]},&f(x)\leq c_\ell\\
\end{array}\right.
\end{equation}
\end{quote}
Let ${\omega}^K$ and $\ell$ satisfy the premise of (!). Note that due to $\Delta_{\ell-1}=\Delta_{\ell-1}^*$ we have $u_\ell=u_\ell^*$, $c_\ell=c_\ell^*$, and $v_\ell=v_\ell^*$, and thus also $\Delta_{\ell\lf}^*=\Delta_{\ell\lf}$, $\Delta_{\ell\rg}^*=\Delta_{\ell\rg}$. Let us consider first the case where the actual estimation procedure terminates due to a disagreement at step $\ell$, so that $\cT^K_{\Delta^*_{\ell\lf}\l}(\bar{\omega}^K)\neq \cT^K_{\Delta^*_{\ell\rg}\r}(\bar{\omega}^K)$. Assuming for a moment that $f(x)< u_\ell=u_\ell^*$, the relation $\bar{\omega}^K\not\in \cE_\ell[x]$ combines with (\ref{badevent}) to imply that $\cT^K_{\Delta^*_{\ell\rg}\r}(\bar{\omega}^K)=\cT^K_{\Delta^*_{\ell\lf}\l}(\bar{\omega}^K)=\dw$, which is impossible under disagreement. Assuming $f(x)> v_\ell=v_\ell^*$, the same argument results in
$\cT^K_{\Delta^*_{\ell\rg}\r}(\bar{\omega}^K)=\cT^K_{\Delta^*_{\ell\lf}\l}(\bar{\omega}^K)=\up$, which again is impossible. We conclude that in the case in question $u_\ell\leq f(x)\leq v_\ell$, i.e., $f(x)\in\bar{\Delta}$, as claimed.\par
Now, assume that there is a consensus at the step $\ell$ in the actual Bisection. When $\bar{\omega}^K\not\in\cE_\ell[x]$ this is only possible when
\begin{enumerate}
\item $\cT^K_{\Delta^*_{\ell\rg}\r}(\bar{\omega}^K)=\dw$ when $f(x)\leq u_\ell=u_\ell^*$,
\item $\cT^K_{\Delta^*_{\ell\rg}\r}(\bar{\omega}^K)=\dw$ when $u_\ell< f(x)\leq c_\ell=c_\ell^*$,
\item $\cT^K_{\Delta^*_{\ell\lf}\l}(\bar{\omega}^K)=\up$ when $c_\ell < f(x)< v_\ell=v_\ell^*$,
\item $\cT^K_{\Delta^*_{\ell\lf}\l}(\bar{\omega}^K)=\up$ when $v_\ell \leq f(x)$,
\end{enumerate}
In situations 1 and 2, and due to consensus at the step $\ell$, (\ref{NDeltaell}) means that $\Delta_\ell=[a_{\ell-1},c_\ell]$, which combines with (\ref{NDeltaellstar}) and $v_\ell=v_\ell^*$ to imply that $\Delta_\ell=\Delta_\ell^*$. Similarly, in situations 3-4 and due to consensus at the step $\ell$, (\ref{NDeltaell}) says that $\Delta_\ell=[c_\ell,b_{\ell-1}]$, which combines with $u_\ell=u_\ell^*$ and (\ref{NDeltaellstar}) to imply that  $\Delta_\ell=\Delta_\ell^*$.
\qed

\subsubsection{Proof of Proposition \ref{propBisection}(ii)}
There is nothing to prove when ${b_0-a_0\over 2}\leq \rho$, since in this case the estimate ${a_0+b_0\over2}$ which does not use observations at all is $(\rho,0)$-reliable. From now on we assume that $b_0-a_0> 2\rho$, implying that $L$ is positive integer.
\paragraph{1$^0$.} Observe, first, that if $a,b$ are such that $a$ is lower-feasible,  $b$ is upper-feasible, and $b-a>2\bar\rho$,  then for every $i\leq I_{b,\geq}$ and $j\leq I_{a,\leq}$ there exists a test, based on $\bar{K}$ observations, which decides upon the hypotheses $H_1$, $H_2$, stating that the observations are drawn from $p_{A(x)}$ with $x\in Z^{b,\geq}_i$ ($H_1$) and
with $x\in Z^{a,\leq}_j$ ($H_2$) with risk at most $\epsilon$. Indeed, it suffices to consider the test which accepts $H_1$ and rejects $H_2$ when $\widehat{f}(\omega^{\bar{K}})\geq{a+b\over 2}$ and accepts $H_2$ and rejects $H_1$ otherwise.
\paragraph{2$^0$.} With parameters of Bisection chosen according to (\ref{Nsetup}), by {Lemma \ref{lem:appendix000}}  we have
\begin{quote}
{(E.1)} {\sl For every $x\in X$, the $p_{A(x)}$-probability of the event $f(x)\in\bar{\Delta}$, $\bar{\Delta}$ being the output segment  of our Bisection, is at least $1-\epsilon$.}
\end{quote}
\paragraph{3$^0$.} We claim that
\begin{enumerate}
\item[(F.1)] {\sl Every segment $\Delta=[a,b]$ with $b-a>2\bar\rho$ and lower-feasible $a$ is $\delta$-good (right)},
\item[(F.2)] {\sl Every segment $\Delta=[a,b]$ with $b-a>2\bar\rho$ and upper-feasible $b$ is $\delta$-good (left)},
\item[(F.3)] {\sl Every $\varkappa$-maximal $\delta$-good (left or right) segment has length at most $2\bar\rho+\varkappa=\rho$. As a result, for every {constructive} step $\ell$, the lengths of the segments $\Delta_{\ell\rg}$ and $\Delta_{\ell\lf}$ do not exceed $\rho$.}
\end{enumerate}
Let us verify (F.1) (verification of F.2 is completely similar, and (F.3) is an immediate consequence of (F.1) and (F.2)). Let $[a,b]$ satisfy the premise of (F.1). It may happen that $b$ is upper-infeasible,  whence $\Delta=[a,b]$ is $0$-good (right), and we are done. Now let $b$ be upper-feasible. As we have already seen, whenever $i\leq I_{b,\geq}$ and $j\leq I_{a,\leq}$, the hypotheses stating that $\omega_k\sim p_{A(x)}$ for some $x\in Z^{b,\geq}_i$, resp., for some $x\in Z^{a,\leq}_j$, can be decided upon with risk $\leq\epsilon$, {implying by (\ref{epsilonstar}) that}
$$
\epsilon_{ij\Delta}\leq [2\sqrt{\epsilon(1-\epsilon)}]^{1/\bar{K}}.
$$
Hence, taking into account that the column and the row sizes of $E_{\Delta\r}$ do not exceed $NI$,
$$
\sigma_{\Delta\r}\leq NI\max_{i,j}\epsilon_{ij\Delta}^K\leq NI[2\sqrt{\epsilon(1-\epsilon)}]^{K/\bar{K}}\leq {\epsilon\over 2L}=\delta
$$
(we have used (\ref{Nsetup})), So, $\Delta$ indeed is $\delta$-good (right).
\paragraph{4$^0$.}
Let us fix $x\in X$ and consider a trajectory of Bisection, the $K$-repeated observation $\omega^K$ being drawn from $p^{K}_{A(x)}$. The output $\bar{\Delta}$ of the procedure is given by one of the following options:
\begin{enumerate}
\item At some step $\ell$ of Bisection, the process terminated {by} \ref{step1r} or \ref{step1l}. In the first case, the segment $[c_\ell,b_{\ell-1}]$ has lower-feasible left endpoint and is not $\delta$-good (right), implying by F.1 that the length of this segment (which is $1/2$ of the length of $\bar{\Delta}=\Delta_{\ell-1}$) is $\leq 2\bar\rho$, so that the length $|\bar{\Delta}|$ of $\bar{\Delta}$ is at most $4\bar\rho\leq 2\rho$. By completely similar argument, the same conclusion holds true when the process terminated at step $\ell$ {by} \ref{step1l}.
\item At some step $\ell$ of Bisection, the process terminated due to disagreement. In this case, by (F.3), we have $|\bar{\Delta}|\leq 2\rho$.
\item Bisection terminated at step $L$, and $\bar{\Delta}=\Delta_L$. In this case, termination clauses {of} \ref{step1r}, \ref{step1l} and \ref{mainstep} were never invoked, clearly implying that $|\Delta_s|\leq {\half}|\Delta_{s-1}|$, $1\leq s\leq L$, and thus $|\bar{\Delta}|=|\Delta_L|\leq {1\over 2^L}|\Delta_0|\leq 2\rho$ (see (\ref{Nsetup})).
\end{enumerate}
Thus, along with (E.1) we have
\begin{quote}
{(E.2)} {\sl It always holds $|\bar{\Delta}|\leq 2\rho$,}
\end{quote}
implying that whenever the signal $x\in X$ underlying observations and the output segment $\bar{\Delta}$ are such that $f(x)\in\bar{\Delta}$, the error of the Bisection estimate (which is the midpoint of $\bar{\Delta}$) is at most $\rho$. Invoking (E.1), we conclude that the Bisection estimate is $(\rho,\epsilon)$-reliable. \qed

\section{1-convexity of conditional quantile}\label{1convex}
{
Let $r$ be a nonvanishing probability distribution on $S$, and let
$$
F_m(r)=\sum_{i=1}^m r_i,\,1\leq m\leq M,
$$
so that $0<F_1(r)<F_2(r)<...<F_M(r)=1$. Denoting by $\cP$ the set of all nonvanishing probability distributions on $S$, observe that for every $p\in\cP$  $\chi_\alpha[r]$ is a piecewise linear function of $\alpha\in [0,1]$ with breakpoints $0,F_1(r),F_2(r),F_3(r),...,F_M(r)$, the values of the function at these breakpoints being $s_1,s_1,s_2,s_3,...,s_M$. In particular, this function is equal to $s_1$ on $[0,F_1(r)]$ and is strictly increasing on $[F_1(r),1]$. Now let $s\in \bR$, and let
$$\cP_\alpha^\leq[s]=\{r\in\cP:\chi_\alpha[r]\leq s\},\;\;\cP_\alpha^\geq[s]=\{r\in\cP: \chi_\alpha[r]\geq s\}.$$
Observe that the just introduced sets are cut off $\cP$ by nonstrict  linear inequalities, specifically,
\begin{itemize}
\item when $s<s_1$, we have $\cP_\alpha^{\leq}[s]=\emptyset$, $\cP_\alpha^{\geq}[s]=\cP$;
\item when $s=s_1$, we have $\cP_\alpha^{\leq}[s]=\{r\in \cP:F_1(r)\geq\alpha\}$, $\cP_\alpha^{\geq}[s]=\cP$;
\item when $s>s_M$, we have $\cP_\alpha^{\leq}[s]=\cP$, $\cP_\alpha^{\geq}[s]=\emptyset$;
\item when $s_1<s\leq s_M$, for every $r\in\cP$ the equation $\chi_\gamma[r]=s$ in variable $\gamma\in[0,1]$ has exactly one solution $\gamma(r)$ which can be
found as follows: we specify $k=k^s\in\{1,...,M-1\}$ such that $s_k<s\leq s_{k+1}$ and set
$$
\gamma(r)={(s_{k+1}-s)F_k(r)+(s-s_k)F_{k+1}(r)\over s_{k+1}-s_k}.
$$
Since $\chi_\alpha[r]$ is strictly increasing in $\alpha$ when $\alpha\in[F_1(p),1]$, for $s\in(s_1,s_M]$ we have
\bse
\cP_\alpha^\leq[s]=\{r\in\cP:\alpha\leq\gamma(r)\}=\left\{r\in\cP:{(s_{k+1}-s)F_k(r)+(s-s_k)F_{k+1}(r)\over s_{k+1}-s_k}\geq\alpha\right\},\\
\cP_\alpha^\geq[s]=\{r\in\cP:\alpha\geq\gamma(r)\}=\left\{r\in\cP:{(s_{k+1}-s)F_k(r)+(s-s_k)F_{k+1}(r)\over s_{k+1}-s_k}\leq\alpha\right\}.
\ese
\end{itemize}
Now, given $\tau\in T$ and $\alpha\in[0,1]$, let us set
\[
G_{\tau,\mu}(p)=\sum_{\iota=1}^\mu p(\iota,\tau),\,1\leq\mu\leq M,
\]
and
\[
\cX^{s,\leq}=\{p(\cdot,\cdot)\in\cX:\chi_\alpha[p_\tau]\leq s\},\;\;\cX^{s,\geq}=\{p(\cdot,\cdot)\in\cX:\chi_\alpha[p_\tau]\geq s\}.
\]
As an immediate consequence of the above description
we get
\bse
s<s_1&\Rightarrow&\cX^{s,\leq}=\emptyset,\;\;\cX^{s,\geq}=\cX,\\
s=s_1&\Rightarrow&\cX^{s,\leq}=\{p\in\cX:G_{\tau,1}(p)\leq s_1G_{\tau,M}(p)\},\;\;\cX^{s,\geq}=\cX,\\
s>s_M&\Rightarrow&\cX^{s,\leq}=\cX,\;\;\cX^{s,\geq}=\emptyset,\\
s_1<s\leq s_M&\Rightarrow&\left\{\begin{array}{l}\cX^{s,\leq}=\left\{p\in\cX:{(s_{k+1}-s)G_{\tau,k}(r)+(s-s_k)G_{\tau,k+1}(r)\over s_{k+1}-s_k}\geq  \alpha G_{\tau,M}(p)\right\},\\
\cX^{s,\geq}=\left\{p\in\cX:{(s_{k+1}-s)G_{\tau,k}(r)+(s-s_k)G_{\tau,k+1}(r)\over s_{k+1}-s_k}\leq \alpha G_{\tau,M}(p)\right\},\\
\end{array}\right.\\
&&k=k_s: s_k<s\leq s_{k+1},
\ese
implying 1-convexity of the conditional quantile on $\cX$ (recall that $G_{\tau,\mu}(p)$ are linear in $p$).
}

\end{document}